\documentstyle[12pt]{article}
\include{amssymb}
\newtheorem{prop}{Proposition}
\newtheorem{def-prop}{Definition-Proposition}
\newtheorem{thm}{Theorem}
\newtheorem{cor}{Corollary}

\newtheorem{lemma}{Lemma}
\newtheorem{defn}{Definition}
\def\endproof{$\Box$}
\def\1{1}
\def\P{{\bf P}}

\def\Z{{\bf Z}}
\def\A{{\bf A}}
\def\G{{\bf G}}

\def\Q{{\bf Q}}

\def\Mgbar{\overline{{\cal M}}_g}
\def\Mbar{\overline{{\cal M}}}

\begin{document}
\title{Equivariant Intersection Theory}
\author{Dan Edidin\thanks{Department of Mathematics, University of Missouri,
Columbia MO 65211. Partially supported by the NSF and UM Research
Board} \mbox{ }and William Graham\thanks{School of Mathematics,
Institute for Advanced Study, Princeton, NJ 08540.
Partially supported by an NSF post-doctoral fellowship}}
\date{}
\maketitle

\section{Introduction}
The purpose of this paper is to develop an equivariant intersection
theory for actions of linear algebraic groups on
schemes and more generally algebraic spaces.
The theory is based on our construction of equivariant Chow
groups. These are algebraic analogues of equivariant cohomology groups
which have all the functorial properties of ordinary Chow groups.  In
addition, they enjoy many of the properties of equivariant cohomology.

Previous work (\cite{Br}, \cite{Gi}, \cite{Vi}) defined equivariant
Chow groups using only invariant cycles on $X$.  However, there are
not enough invariant cycles on $X$ to define equivariant Chow groups
with nice properties, such as being a homotopy invariant, or having an
intersection product when $X$ is smooth (see Section \ref{noinprod}).
The definition of this paper is modeled after Borel's
definition of equivariant cohomology. It is made possible
by Totaro's approximation of $EG$ by open subsets
of representations (\cite{To}). Consequently,  an equivariant class is
represented by an invariant cycle on $X \times V$, where $V$ is a
representation of $G$.  By enlarging the definition of equivariant
cycle, we obtain a rich theory, which is closely related to other
aspects of group actions on schemes and algebraic spaces.

After establishing the basic properties of equivariant Chow groups,
this paper is mainly devoted to the relationship between equivariant
Chow groups and Chow groups of quotient algebraic spaces and stacks.
If $G$ is a linear algebraic group acting on a space $X$, denote by
$A^G_{i}(X)$ the $i$-th equivariant Chow group of $X$.  If $G$ acts
properly on $X$ then a quotient $X/G$ exists in the
category of algebraic spaces (\cite{Kollar}, \cite{KM}); under some
additional hypotheses (see \cite{GIT}) if $X$ is a
scheme then $X/G$ is a scheme.  We prove that there is an
isomorphism of $A^G_{i + \mbox{\small{dim}}G}(X) \otimes \Q$ and $A_i(X/G)
\otimes \Q$.  If $G$ acts with trivial stabilizers this holds without
tensoring with $\Q$.

For an action which is not proper, there need not be a quotient in the
category of algebraic spaces. However, there is always an
Artin quotient stack $[X/G]$.  We prove that the
equivariant groups $A^G_*(X)$ depend only on the stack $[X/G]$ and not
on its presentation as a quotient.  If $X$ is smooth, then $A^1_G(X)$
is isomorphic to Mumford's Picard group of the stack, and the ring
$A^*_G(X)$ can naturally be identified as an integral Chow ring of
$[X/G]$ (Section \ref{intstack}).

These results imply that equivariant Chow groups are a useful tool for
computing Chow groups of quotient spaces and stacks.  For example,
Pandharipande (\cite{Pa1}, \cite{Pa2}) has used equivariant methods to
compute Chow rings of moduli spaces of maps of projective spaces as
well as the Hilbert scheme of rational normal curves.  In this paper,
we compute the integral Chow rings of the stacks ${\cal M}_{1,1}$ and
$\Mbar_{1,1}$ of elliptic curves, and obtain a simple proof of
Mumford's result (\cite{Mu}) that $Pic_{fun}({\cal M}_{1,1}) =
\Z/12\Z$.  In an appendix to this paper, Angelo Vistoli computes the
Chow ring of ${\cal M}_2$, the moduli stack of smooth curves of genus
2.

Equivariant Chow groups are also useful in proving results about
intersection theory on quotients.  It is easy to show that if $X$ is
smooth then there is an intersection product on $A^G_*(X)$.  The
theorem on quotients therefore implies that there exists an
intersection product on the rational Chow groups of a quotient of a
smooth algebraic space by a proper action.  The existence of such an
intersection product was shown by Gillet and Vistoli, but only under
the assumption that the stabilizers are reduced.  This is automatic in
characteristic $0$, but typically fails in characteristic $p$.  The
equivariant approach does not require this assumption and therefore
extends the work of Gillet and Vistoli to arbitrary characteristic.
Furthermore, by avoiding the use of stacks, the proof becomes much
simpler.

Finally, equivariant Chow groups define invariants of quotient stacks
which exist in arbitrary degree, and associate to a smooth quotient
stack an integral intersection ring which when tensored with $\Q$
agrees with rings defined by Gillet and Vistoli.  By analogy with
quotient stacks, this suggests that there should be an integer
intersection ring associated to an arbitrary smooth stack, which could
be nonzero in degrees higher than the dimension of the stack.

We remark that besides the properties mentioned above, the equivariant
Chow groups we define are compatible with other equivariant theories
such as cohomology and $K$-theory.  For instance, if $X$ is smooth
then there is a cycle map from equivariant Chow theory to equivariant
cohomology (Section \ref{cycles}).  In addition, there is a map from
equivariant $K$-theory to equivariant Chow groups, which is an
isomorphism after completion; and there is a localization theorem for
torus actions, which can be used to give an intersection theoretic
proof of residue formulas of Bott and Kalkman.  These topics will be
treated elsewhere.

\tableofcontents

\medskip

{\bf Acknowledgments:} We thank William Fulton,
Rahul Pandharipande and Angelo Vistoli for advice and encouragement.
We also benefited from discussions with Burt Totaro, Amnon Yekutieli,
Robert Laterveer and Ruth Edidin.
Thanks to Holger Kley for suggesting the inclusion
of the cycle map to equivariant cohomology, and to J\'anos
Koll\'ar for emphasizing the algebraic space point of view.

\section{Definitions and basic properties}
\subsection{Conventions and Notation}

This paper is written in the language of algebraic spaces.  It is
possible to work entirely in the category of schemes, provided the
mild technical hypotheses of Proposition \ref{inap} are satisfied.
These hypotheses insure that if $X$ is a scheme with a $G$-action, then the
mixed spaces $X_G = (X \times U)/G$ are schemes.  Here $G$ acts freely
on $U$ (which is an open subset of a representation of $G$) and hence
on $X \times U$.  In the category of algebraic spaces, quotients
of free actions always exist (\cite{D-M} or \cite{Artin}; see
Proposition \ref{l.algspacequotient}).  By working with algebraic
spaces we can therefore avoid the hypotheses of Proposition \ref{inap}.

Another reason to work with algebraic spaces comes from a theorem of
Koll\'ar and Keel and Mori (\cite{Kollar}, \cite{KM}), generalizing
the result about free actions mentioned above.  They prove that if
$G$ is a linear algebraic group acting properly on a separated
algebraic space $X$, then a geometric quotient $X/G$ exists as a
separated algebraic space.  Such quotients arise frequently in moduli
problems.  However, even if $X$ is a scheme, the quotient $X/G$ need
not be a scheme.  By developing the theory for algebraic spaces, we
can apply equivariant methods to study such quotients.

For these reasons, the natural category for this theory is that of
algebraic spaces.  In Section \ref{algspace} we explain why
the intersection theory of \cite{Fulton} remains unchanged in
this category.

Except in Section \ref{mixed}, all schemes and algebraic spaces are
assumed to be quasi-separated and of finite type over a field $k$
which can have arbitrary characteristic.  A {\em smooth} space is
assumed to be separated, implying that the diagonal $X \rightarrow X
\times X$ is a regular embedding. For brevity, we will use the term
variety to mean integral algebraic space (rather than integral {\em
scheme}, as is usual).  An algebraic group is always assumed to be
linear.  For simplicity of exposition, we will usually assume that our
spaces are equidimensional.

\paragraph{Group actions}
If an algebraic group $G$ acts on a scheme or algebraic space $X$ then
the action is said to be {\it closed} if the orbits of geometric
points are closed in $X$.  It is {\it proper} if the action map $G
\times X \rightarrow X \times X$ is proper. If every
point has an invariant neighborhood such that the action
is proper in the neighborhood then we say the action is
{\em locally proper}.
It is
{\it free} if the action map is a closed embedding.

If $G \times X \stackrel{j} \rightarrow X \times X$ is a group action, we
will call $j^{-1}(\Delta_X) \rightarrow \Delta_X$ the stabilizer
group scheme of $X$. Its fibers are the stabilizers of the points
of $X$. A group action is said to have {\em
finite stabilizer} if
the map $j^{-1}(\Delta_X) \rightarrow \Delta_X$ is finite.
Since a linear algebraic group is affine, the
fibers of $j: G \times X \stackrel{j} \rightarrow X \times X$ are
finite when $G$ acts properly. Hence, a proper action
has finite stabilizer. The converse need not be true (\cite[Example 0.4]{GIT}).
Note that a locally proper action has finite stabilizer as well.

Finally we say the action is {\em set theoretically free} or
{\em has trivial stabilizer}
if the stabilizer of every point is trivial. In particular
this means $j^{-1}(\Delta_X) \rightarrow \Delta_X$ is an isomorphism.
If the action is proper and has trivial stabilizer then it is
free (Lemma \ref{l.free}).

A flat, surjective, equivariant map $X \stackrel{f}\rightarrow Y$ is
called a principal bundle if $G$ acts trivially on $Y$, and the map $X
\times G \rightarrow X \times_Y X$ is an isomorphism.  This condition
is equivalent to local triviality in the \'etale topology, i.e., there
is an \'etale cover $U \rightarrow Y$ such that $U \times_Y X \simeq U
\times G$.

As noted above, if $G$ acts set-theoretically freely on $X$ then a
geometric quotient $X/G$ always exists in
the category of algebraic spaces, and moreover, $X$ is a principal
$G$-bundle over $X/G$ (Proposition \ref{l.algspacequotient}).

\subsection{Equivariant Chow groups} \label{basicdef}
Let $X$ be an $n$-dimensional
algebraic space.
We will denote the $i$-th equivariant Chow group
of $X$ by $A^G_i(X)$, and define it as follows.

Let $G$ be a $g$-dimensional algebraic group.  Choose an
$l$-dimensional representation $V$ of $G$ such that $V$ has an open
set $U$ on which $G$ acts freely and whose complement has codimension
more than $n-i$.  Let $U \rightarrow U/G$ be the principal bundle
quotient.  Such a quotient automatically exists as an algebraic space;
moreover, for any algebraic group, representations exist so that $U/G$
is a scheme -- see Lemma \ref{q.exist} of Section \ref{appendix}.  The
principal bundle $U \rightarrow U/G$ is Totaro's finite dimensional
approximation of the classifying bundle $EG \rightarrow BG$ (see
\cite{To} and \cite{E-G}).  The diagonal action on $X \times U$ is
also free, so there is a quotient in the category of algebraic spaces
$X \times U \rightarrow (X \times U)/G$ which is a principal $G$
bundle.  We will usually denote this quotient by $X_G$. (See
Proposition \ref{inap} for conditions that are sufficient to imply
that $X_G$ is a scheme.)

\begin{def-prop} \label{keydef}
Set $A_i^G(X)$ (the $i$-th equivariant Chow group) to be
$A_{i+l-g}(X_G)$, where $A_*$ is the usual Chow group.
This group is independent of the representation
as long as $V- U$ has sufficiently high codimension.
\end{def-prop}

\medskip

Proof.  As in \cite{To}, we
will use Bogomolov's double fibration argument.  Let $V_1$ be another
representation of dimension $k$ such that there is an open subset
$U_1$ with a principal bundle quotient $U_1 \rightarrow U_1/G$ and
whose complement has codimension at least $n-i$.  Let $G$ act
diagonally on $V \oplus V_1$.  Then $V \oplus V_1$ contains an open
set $W$ which has a principal bundle quotient $W/G$ and contains both
$U \oplus V_1$ and $V \oplus U_1$.  Thus, $A_{i+k+l-g}(X \times^G W) =
A_{i+k+l-g}(X \times^G (U \oplus V_1))$ since $(X \times^G W)-(X
\times^G (U \oplus V_1)$ has dimension smaller than $i+k+l-g$. On the
other hand, the projection $V \oplus V_1 \rightarrow V$ makes $X
\times^G (U \oplus V_1)$ a vector bundle over $X \times^G U$ with
fiber $V_1$ and structure group $G$. Thus, $A_{i+k+l-g}(X \times^G (U
\oplus V_1)) = A_{i+l-g}(X \times^G U)$. Likewise, $A_{i+k+l-g}(X
\times^G W) =A_{i+k-g}(X \times^G U_1)$, as desired.  \endproof \medskip

{\bf Remark.} In the sequel, the notation $U \subset V$ will refer to
an open set in a representation on which the action is free, and
$X_G$ will mean a mixed quotient $X \times^G U$ for any
representation $V$ of $G$. If we write $A_{i+l-g}(X_G)$ then $V-U$ is
assumed to have codimension more than $n-i$ in $V$. (As above
$n=\mbox{dim }X$, $l=\mbox{dim }V$ and $g =\mbox{dim }G$.)

\paragraph{Equivariant cycles}
If $Y \subset X$ is an $m$-dimensional $G$-invariant subvariety
(recall that variety means integral algebraic space), then
it has a $G$-equivariant fundamental class $[Y]_G \in A_m^G(X)$.
More generally, if $V$ is an $l$-dimensional representation
and $S \subset X \times V$ is an $m+l$-dimensional subvariety,
then $S$ has a $G$-equivariant fundamental class $[S]_G \in A_m^G(X)$.
Thus, unlike ordinary Chow groups, $A_i^G(X)$
can be non-zero for any $i \leq n$, including negative $i$.

\begin{prop} \label{eqcycle}
If $\alpha \in A_m^G(X)$, then there exists a representation
$V$ such that $\alpha = \sum a_i[S_i]_G$, where $S_i$ are $m +l$
invariant subvarieties of $X \times V$, where $l$ is the dimension of
$V$.
\end{prop}

Proof. Cycles of dimension $m+l-g$ on $X_G$ correspond exactly to
invariant cycles of dimension $m+l$ on $X \times U$. Since $V-U$ has
high codimension, invariant $m+l$ cycles on $X \times U$ extend
uniquely to invariant $m+l$ cycles on $X \times V$. \endproof \medskip

The representation $V$ is not unique. For example, $[X]_G$ and $[X
\times V]_G$ define the same equivariant class.

The projection $X \times U \rightarrow U$ induces
a flat map $X_G \rightarrow U$ with fiber $X$. Restriction
to a fiber gives a map $i^*:A_*^G(X) \rightarrow A_*(X)$
from equivariant Chow groups to ordinary Chow groups. The
map is independent of the choice of fiber because any two
points of $U/G$ are rationally equivalent.
For any $G$-invariant subvariety $Y \subset X$,
$i^*([Y]_G) =[Y]$.

Before reading further, the reader may want to skip to Section \ref{examples}
for examples.

\subsection{Functorial properties}
In this section all maps $f: X \rightarrow Y$ are assumed to
be $G$-equivariant.

Let $\P$ be one of the following properties of morphisms
of schemes: proper, flat, smooth, regular embedding or l.c.i.

\begin{prop}
If $f: X \rightarrow Y$ has property $\P$, then
the induced map $f_G: X_G \rightarrow Y_G$
also has property $\P$.
\end{prop}
Proof. If $X \rightarrow Y$ has property $\P$, then, by base change,
so does the map $X \times U \rightarrow Y \times U$.
The morphism $Y \times U \rightarrow Y_G$ is flat and surjective
(hence faithfully flat), and
$X \times U \simeq X_G \times_{Y_G} Y \times U$. Thus by descent
\cite[Section 8.4 - 5]{SGA1}, the morphism $X_G \rightarrow Y_G$
also has property $\P$.
\endproof

\begin{prop}
Equivariant Chow groups have the same functoriality
as ordinary Chow groups for equivariant morphisms with property
$\P$.
\end{prop}
Proof. If $f:X \rightarrow Y$ has property $\P$, then so
does $f_G:X_G \rightarrow Y_G$. Define pushforward $f_*$ or pullback
$f^*$ on equivariant Chow groups as the pullback or pushforward
on the ordinary Chow groups of $X_G$ and $Y_G$. The double fibration
argument shows that this is independent of choice of representation.
\endproof

\subsection{Chern classes}
Let $X$ be a scheme with a $G$-action, and let
$E$ be an equivariant vector bundle (in the category of algebraic
spaces)
Consider the quotient $E \times U \rightarrow E_G$.

\begin{lemma}
$E_G \rightarrow X_G$ is a vector bundle.
\end{lemma}
Proof.  The bundle $E_G \rightarrow X_G$ is an affine bundle which is
locally trivial in the \'etale topology since it becomes locally
trivial after the smooth base change $X \times U \rightarrow
X_G$. Also, the transition functions are linear since they are linear
when pulled back to $X \times U$.  Hence, $E_G \rightarrow X_G$ is a
vector bundle over $X_G$.  \endproof

\begin{defn}
Define equivariant Chern classes $c_j^G(E):A_i^G(X) \rightarrow A_{i-j}^G(X)$
by $c_j^G(E)\cap \alpha= c_j(E_G) \cap \alpha \in A_{i-j+l-g}(X_G)$.
\end{defn}
By the double fibration argument, the definition does not depend on the choice
of representation.

Following \cite{GIT}, we denote by $Pic^G(X)$ the group of
isomorphism classes of $G$-linearized locally free sheaves on $X$.

\begin{thm} \label{piciscool}
Let $X$ be a locally factorial variety of dimension $n$.
Then the map $\mbox{Pic}^G(X) \rightarrow A_{n-1}^G(X)$
defined by $L \mapsto (c_1(L) \cap [X]_G)$ is an isomorphism.
\end{thm}
Proof.
We know that the map $Pic(X_G) \stackrel{\cap c_1(L_G)}
\rightarrow A_{n-g+l-g}(X_G)= A_{n-1}^G(X)$ is an isomorphism. Since
$X \times U \rightarrow X_G$ is a principal bundle,
$Pic(X_G) = Pic^G(X \times U)$. The theorem now follows
from the following lemma.

\begin{lemma} \label{units}
Let $X$ be a locally factorial variety with a $G$-action.

(a) Let $U \stackrel{j} \hookrightarrow X$ be an invariant subvariety
such that $X-U$ has codimension more than 1.
Then the restriction map
$j^*:Pic^G(X) \rightarrow Pic^G(U)$ is an isomorphism.

(b) Let $V$ be a representation and let $\pi: X \times V \rightarrow
X$ be the projection.  Then
$\pi^*:Pic^G(X) \rightarrow Pic^G(X \times V)$ is an
isomorphism.
\end{lemma}
Proof of Lemma \ref{units}.
We first prove (a).

Injectivity: Suppose $L \in \mbox{Pic}^G(X)$ and $j^*L$ is
trivial.  Since $\mbox{Pic}(X) \cong \mbox{Pic}(X-Y)$, this implies that
as a bundle $L$ must be trivial.  A linearization of the trivial bundle
on $X$ is just a homomorphism $G \rightarrow \Gamma(X,{\cal O}^*_X)$.
Since $X$ is a variety and
$X-U$ has codimension more than one, $\Gamma(X,{\cal O}^*_X) =
\Gamma(U,{\cal O}^*_U)$. Thus a
linearization of the trivial bundle is trivial on $X$ if and only if
it is trivial on $U$, proving injectivity.

Surjectivity: A linearization of $L$ is a homomorphism of $G$ into
the group of automorphisms of $L$ over $X$.
To show that $j^*$ is surjective, we must
show that if $L \in \mbox{Pic}(X)$ is linearizable on $U$ then it is
linearizable on $X$.  But any isomorphism $\alpha: L |_{U} \rightarrow
g^*L|_{U}$ extends to an isomorphism over $X$.  (To see this, pick an
isomorphism $\beta: L \rightarrow g^*L$ ; we know one exists because
$\mbox{Pic}(X) \cong \mbox{Pic}(U)$ and $L$ and $g^*L$ are isomorphic
on $X-Y$.  Then $\alpha = \beta \cdot f$, where $f \in  \Gamma(U,{\cal
O}^*(U))$, but $\Gamma(X,{\cal O}^*_X) = \Gamma(U, {\cal O}^*_U)$,
so $\alpha$
extends to $X$.)  Hence $L$ is linearizable on $X$.

The proof of (b) is similar. The key point is that if
$X$ is a variety and $V$ is a vector space, then $\Gamma(X \times V,
{\cal O}^*_{X \times V}) = \Gamma(X, {\cal O}^*_X)$, because if
$R$ is an integral domain, then the units in $R[t_1, \ldots , t_n]$
are the just the units of $R$.
\endproof

\subsection{Exterior Products}
If $X$ and $Y$ have $G$-actions then there are exterior
products $A_i^G(X) \otimes A_j^G(X) \rightarrow A^G_{i+j}(X \times Y)$.
By Proposition \ref{eqcycle} any $\alpha \in A_*^G(X)$ can
be written as $\alpha = \sum a_i [S_i]_G$ where the $S_i$'s
are $G$-invariant subvarieties of $X \times V$ for some representation
$V$.

Let $V$, $W$ be representations
of $G$ of dimensions $l$ and $k$ respectively. Let
$S \subset X \times V$, $T \subset Y \times W$, be $G$-invariant subvarieties
of dimensions $i+l$ and $j +k$ respectively. Let $s: X \times V \times Y \times
W \rightarrow X \times Y \times (V \oplus W)$ be the isomorphism
$(x,v,y,w) \mapsto (x,y,v \oplus w)$.
\begin{def-prop} (Exterior products)
The assignment $[S]_G \times [T]_G \mapsto [s(S \times T)]_G$
induces a well defined exterior product map of equivariant Chow groups
$A_i^G(X) \otimes A_j^G(Y) \rightarrow A_{i+j}^G(X \times Y)$.
\end{def-prop}
Proof. The proof follows from \cite[Proposition 1.10]{Fulton} and
the double fibration argument used above. \endproof \medskip

Given the above propositions, equivariant Chow groups
satisfy all the formal properties of ordinary Chow groups
(\cite[Chapters 1-6]{Fulton}). In particular, if $X$ is
smooth, there is an intersection product on the
the equivariant Chow groups $A_*^G(X)$ which makes $\oplus A_*^G(X)$
into a graded ring.

\subsection{Operational Chow groups}
In this section we define equivariant operational Chow groups
$A^i_G(X)$, and compare them with the operational Chow groups of $X_G$.

Define equivariant operational Chow groups $A^i_G(X)$ as operations
$c(Y \rightarrow X): A_*^G(Y) \rightarrow A_{*-i}^G(Y)$ for every
$G$-map $Y \rightarrow X$.  As for ordinary operational Chow groups
(\cite[Chapter 17]{Fulton}), these operations should be compatible
with the operations on equivariant Chow groups defined above (pullback
for l.c.i. morphisms, proper pushforward, etc.).  From this definition
it is clear that for any $X$, $A^*_G(X)$ has a ring structure.  The
ring $A^*_G(X)$ is positively graded, and $A^i_G(X)$ can be non-zero
for any $i \geq 0$.

Note that by construction, the equivariant Chern classes defined above
are elements of the equivariant operational Chow ring.

\begin{prop} \label{opsmooth}
If $X$ is smooth of dimension $n$,
then $A^i_G(X) \simeq A_{n-i}^G(X)$.
\end{prop}

\begin{cor} (of Theorem \ref{piciscool})
If $X$ is a smooth variety with a $G$-action, then the map
$Pic^G(X) \rightarrow A^1_G(X)$ defined by $L \mapsto c_1(L)$
is an isomorphism. \endproof
\end{cor}

Proof of Proposition \ref{opsmooth}.
Define a map $A^i_G(X) \rightarrow A_{n-i}^G(X)$
by the formula $c \mapsto c \cap [X]_G$.
Define a map
$A_{n-i}^G(X)\rightarrow A^i_G(X)$, $\alpha \mapsto c_{\alpha}$
as follows.  If $Y \stackrel{f} \rightarrow X$ is a
$G$-map, then since $X$ is smooth, the graph $\gamma_f: Y \rightarrow
Y \times X$ is a $G$-map which is a regular embedding.  If $\beta \in
A^G_*(Y)$ set $c_\alpha \cap \beta= \gamma_f^*(\beta \times
\alpha)$.

\medskip

Claim (cf. \cite[Proposition 17.3.1]{Fulton}): $\beta \times
(c \cap [X]_G)
 = c \cap (\beta \times  [X]_G)$.

\medskip
Given the claim, the formal arguments of
\cite[Proposition 17.4.2]{Fulton}
show that the two maps are inverses.

Proof of Claim: By Proposition \ref{eqcycle} and the linearity
of of equivariant operations, we may assume there is a representation $V$
so that  $\beta = [S]_G$ for a $G$-invariant subvariety $S \subset Y \times V$.
Since $S$ is $G$-invariant,
the projection $p: S \times X \rightarrow X$ is equivariant.
Thus,
$$[S]_G \times (c \cap [X]_G) = p^*(c \cap [X]_G)=
c \cap p^*([X]_G) = c \cap ([S \times X]_G) = c \cap ([S]_G \times [X]_G).$$
\endproof

\medskip

Let $V$ be a representation such that
$V- U$ has codimension more than $k$, and set $X_G =
X \times^G U$. In the remainder of the subsection we will
discuss the relation between $A^k_G(X)$ and $A^k(X_G)$ (ordinary
operational Chow groups).

Recall \cite[Definition 18.3]{Fulton} that an envelope
$\pi:\tilde{X} \rightarrow X$ is a proper map such that for any
subvariety $W \subset X$ there is a subvariety $\tilde{W}$ mapping
birationally to $W$ via $\pi$. In the case of group actions, we will
say that $\pi: \tilde{X} \rightarrow X$ is an {\it equivariant}
envelope, if $\pi$ is $G$-equivariant, and if we can take $\tilde{W}$
to be $G$-invariant for $G$-invariant $W$. If there is an open set $X^0
\subset X$ over which $\pi$ is an isomorphism, then we say $\pi:
\tilde{X} \rightarrow X$ is a {\it birational} envelope.

\begin{lemma} If $\pi: \tilde{X} \rightarrow X$ is an
equivariant (birational)  envelope, then
$p: \tilde{X}_G \rightarrow X_G$ is a (birational)
 envelope ($\tilde{X}_G$ and $X_G$ are constructed
with respect to a fixed representation $V$). Furthermore,
if $X^0$ is the open set over which $\pi$ is an isomorphism
(necessarily $G$-invariant), then $p$ is an isomorphism
over $X^0_G = X^0 \times^G U$.
\end{lemma}
Proof. Fulton \cite[Lemma 18.3]{Fulton} proves that
the base extension of an envelope is an envelope.
Thus $\tilde{X} \times U \stackrel{\pi \times id}\rightarrow X \times U$
is an envelope. Since the projection $X \times U \rightarrow X$
is equivariant, this envelope is also equivariant.
If $W \subset X_G$ is a subvariety, let $W'$ be its inverse image
(via the quotient map) in $X \times U$. Let $\tilde{W'}$ be
an invariant subvariety of $\tilde{X} \times U$ mapping
birationally to $W'$. Since $G$ acts freely on $\tilde{X} \times U$
it acts freely on $\tilde{W'}$, and $\tilde{W} = \tilde{W'}/G$
is a subvariety of $\tilde{X}_G$ mapping birationally to $W$.
This shows that $\tilde{X}_G \rightarrow X_G$ is an envelope;
it is clear that the induced map $\tilde{X}_G \rightarrow
\tilde{X}$ is an isomorphism over $X_0^G$. \endproof

\medskip

Suppose $\tilde{X} \stackrel{\pi}\rightarrow X$
is an equivariant envelope which is
an isomorphism over $U$. Let $\{S_i\}$ be the irreducible components
of $S= X -X^0$, and let $E_i = \pi^{-1}(S_i)$. Then $\{S_{i G}\}$
are the irreducible components of $X_G - X^0_G$ and
$E_{i G} = \pi^{-1}(S_{i G})$.

\begin{thm}
If $X$ has an equivariant smooth envelope
$\pi: \tilde{X} \rightarrow X$ such that there is an
open $X^0 \subset X$ over which $\pi$ is an isomorphism,
and $V-U$ has codimension more than $k$, then
$A^k_G(X) = A^k(X_G)$.
\end{thm}

Proof. If $\pi: \tilde{X} \rightarrow X$ is an
equivariant non-singular  envelope, then
$p: \tilde{X}_G \rightarrow X_G$
is also an  envelope and $\tilde{X}_G$ is non-singular.
Thus, by \cite[Lemma 1.2]{Kimura}
$p^*:A^*(X_G) \rightarrow A^*(\tilde{X}_G)$ is injective.
The image of $p^*$ is described inductively
in \cite[Theorem 3.1]{Kimura}. A class
$\tilde{c} \in A^*(\tilde{X}_G)$ equals
$p^*c$ if and only if for each
$E_{i G}$ , $\tilde{c}_{| E_{i G}} = p^*c_i$
where $c_i \in A^*(E_i)$.
This description follows from formal properties of operational
Chow groups, and the exact sequence \cite[Theorem 2.3]{Kimura}

$$A^*(X_G) \stackrel{p}\rightarrow A^*(\tilde{X}_G)
\stackrel{p_1^* - p_2^*} \rightarrow A^*(\tilde{X}_G \times_{X_G}
\tilde{X}_G)$$ where $p_1$ and $p_2$ are the two projections
from $\tilde{X}_G \times_{X_G} \tilde{X}_G$.

By Proposition \ref{opsmooth} above, we know that
$A^k_G(\tilde{X}) = A^k(\tilde{X}_G)$.
We will show that $A^k_G(X)$ and $A^k(X_G)$ have the same image
in $A^k(\tilde{X}_G)$.
By Noetherian induction we may assume that
$A^k_G(S_i) = A^k((S_{i})_G)$. To prove the theorem, it suffices
to show that there is also an exact sequence of equivariant
operational Chow groups
$$0 \rightarrow A^*_G(X) \stackrel{\pi^*}\rightarrow A^*_G(\tilde{X})
\stackrel{p_{1}^* -p_{2}^*}\rightarrow A^*_G(\tilde{X} \times_X
\tilde{X})$$
This can be checked by working with the action of $A^*_G(X)$
on a fixed Chow group $A_{i}(X_G)$ and arguing as in Kimura's
paper.
\endproof

\begin{cor}
If $X$ is separated and has an
equivariant resolution of singularities (in particular
if the characteristic is 0), and $V-U$ has codimension more than $k$,
then
$A^k_G(X) = A^k(X_G).$
\end{cor}
Proof (cf. \cite[Remark 3.2]{Kimura}).
We must show the existence of an equivariant envelope
$\pi:\tilde{X} \rightarrow X$. By equivariant
resolution of singularities, there is a resolution
$\pi_1:\tilde{X_1} \rightarrow X$ such
that $\pi_1$ is an isomorphism outside
some invariant subscheme $S \subset X$. By Noetherian
induction, we may assume that we have constructed an
equivariant envelope $\tilde{S} \rightarrow S$. Now
set $\tilde{X} = \tilde{X_1} \cup \tilde{S}$.
\endproof \medskip

\subsection{Equivariant higher Chow groups} \label{higheq}
In this section assume that $X$ is quasi-projective
(a quasi-projective algebraic space is a scheme (\cite[p.140]{Knutson}).
Bloch (\cite{Bl}) defined
higher Chow groups $A^i(X,m)$ as $H_m(Z^i(X,\cdot))$
where $Z^i(X,\cdot)$ is a complex whose $k$-th term
is the group of cycles of codimension $i$ in $X \times \Delta^k$
which intersect the faces properly. Since we prefer
to think in terms of dimension rather than codimension
we will define $A_p(X,m)$ as $H_m(Z_p(X,\cdot))$,
where $Z_p(X,k)$ is
the group of cycles of dimension $p+k$ in $X \times \Delta^k$
intersecting the faces properly. When $X$ is equidimensional
of dimension $n$, then $A_p(X,m) = A^{n-p}(X,m)$.

If $Y \subset X$ is closed, there is a localization long exact sequence.
The advantage of indexing by dimension rather than codimension is that
the sequence exists without assuming that $Y$ is equidimensional.

\begin{lemma}
Let $X$ be equidimensional, and let $Y \subset X$ be closed,
then there is a long exact sequence of higher Chow groups
$$\ldots \rightarrow A_p(Y,k) \rightarrow A_p(X,k) \rightarrow
A_p(X-Y,k) \rightarrow \\
\ldots \rightarrow A_p(Y) \rightarrow A_p(X) \rightarrow A_p(X-Y)
\rightarrow 0$$
(there is no requirement that $Y$ be equidimensional).
\end{lemma}
Proof. This is a simple consequence of
the localization theorem of \cite{Bl}.
We must show that the complex $Z_p(X - Y,\cdot)$
is quasi-isomorphic to the complex $\frac{Z_p(X,\cdot)}{Z_p(Y,\cdot)}$.
By induction on the number of components, it suffices to verify the
quasi-isomorphism when $Y$
is the union of two irreducible components $Y_1$ and $Y_2$.

By the original localization theorem, $Z_p(X-(Y_1 \cup Y_2),\cdot)
\simeq \frac{Z_p(X-Y_1,\cdot)}{Z_p(Y_2-(Y_1 \cap Y_2),\cdot)}$ and
$Z_p(X-Y_1, \cdot) \simeq \frac{Z_p(X,\cdot)}{Z_p(Y_1)}$ (here
$\simeq$ denotes quasi-isomorphism).  By induction on dimension, we
can assume that the lemma holds for schemes of smaller dimension, so
$Z_p((Y_2 - (Y_1 \cap Y_2),\cdot) \simeq \frac{Z_p(Y_2,\cdot)}{Z_p(Y_1
\cap Y_2)}$. Finally note that $\frac{Z_p(Y_2,\cdot)}{Z_p(Y_1 \cap
Y_2,\cdot)} \simeq \frac{Z_p(Y_1 \cup Y_2,\cdot)}{Z_p(Y_1,\cdot)}$.
Combining all our quasi-isomorphisms we have
$$Z_p(X-(Y_1 \cup Y_2), \cdot) \simeq \frac{\frac{Z_p(X,\cdot)}
{Z_p(Y_1,\cdot)}}{\frac{Z_p(Y_1 \cup Y_2),\cdot}{Z_p(Y_1,\cdot)}}
\simeq \frac{Z_p(X,\cdot)}{Z_p(Y_1 \cup Y_2,\cdot)}$$
as desired.
\endproof

\medskip

If $X$ is quasi-projective with a linearized $G$-action, we can define
equivariant higher Chow groups $A_{i}^G(X,m)$ as $A_{i+l-g}(X_G,m)$,
where $X_G$ is formed from an $l$-dimensional representation $V$ such
that $V-U$ has high codimension (note that $X_G$ is quasi-projective,
by \cite[Prop. 7.1]{GIT}).  The homotopy property of higher Chow
groups shows that $A_{i}^G(X,m)$ is well defined.

{\bf Warning.} Since the homotopy property of higher Chow groups has
only been proved for quasi-projective varieties, our definition of
higher equivariant Chow groups is only valid for quasi-projective
varieties with a linearized action.  However, if $G$ is connected and
$X$ is quasi-projective and normal, then by Sumihiro's
equivariant completion \cite{Sumihiro} and \cite[Corollary 1.6]{GIT},
any action is linearizable.

\medskip

One reason for constructing equivariant higher Chow groups
is to obtain a localization exact sequence:
\begin{prop} Let $X$ be equidimensional and quasi-projective with a
linearized $G$-action, and
let $Y \subset X$ be an invariant subscheme.  There is a long exact
sequence of higher equivariant Chow groups
$$\ldots \rightarrow A_p^G(Y,k) \rightarrow A_p^G(X,k) \rightarrow
A_p^G(X-Y,k) \rightarrow \\
\ldots \rightarrow A_p^G(Y) \rightarrow A_p^G(X) \rightarrow A_p^G(X-Y)
\rightarrow 0. $$
\endproof
\end{prop}

\subsection{Cycle Maps} \label{cycles}
If $X$ is a complex algebraic variety with the action of
a complex algebraic group, then we can define
equivariant Borel-Moore homology $H_{BM, i}^G(X)$
as $H_{BM,i+2l-2g}(X_G)$ for  $X_G = X \times^G U$.
As for Chow groups, the definition is independent
of the representation, as long as $V -U$ has sufficiently
high codimension,  and we obtain a cycle map
$$cl:A^G_i(X) \rightarrow H_{BM,2i}^G(X)$$
compatible with the usual operations on equivariant
Chow groups (cf. \cite[Chapter 19]{Fulton}).

Let $EG \rightarrow BG$ be the classifying bundle.
The open subsets $U \subset V$ are topological
approximations to $EG$. For, if $\phi$ is a map of the $j$-sphere
$S^j$ to $U$, we may view $\phi$ as a map $S^j \rightarrow V$.  Extend
$\phi$ to a map $B^{j+1} \rightarrow V$.  We may assume that the extended map
is smooth and transversal to $V-U$.  If $j+1 < 2i$, where $i$ is the
complex codimension of $V - U$, then transversality implies that the
extended map does not intersect $V-U$.  Thus we have extended $\phi$
to a map $B^{j+1} \rightarrow U$.  Hence $\pi_j(U) = 0$ for $j <
2i-1$.

Note that $H_{BM,i}^G(X)$ is not the same as $H_i(X \times^G EG)$,
However, if $X$ is smooth, then $X_G$ is also smooth, and $H_{BM,i}(X_G)$
is dual to $H^{2n-i}(X_G)=H^{2n-i}(X \times^G EG)=H^{2n-i}_G(X)$,
where $n$ is the complex dimension of $X$. In this
case we can interpret the cycle
map as giving a map
$$cl: A^i_G(X) \rightarrow H^{2i}_G(X).$$

If $X$ is compact, and the open sets $U \subset V$ can be chosen so
that $U/G$ is projective, then
Borel-Moore homology of $X_G$ coincides with ordinary
homology, so $H^G_{BM*}(X)$ can be calculated with a compact model.
In general, however, $U/G$ is only quasi-projective.
If $G$ is finite, then $U/G$ is
never projective.  If $G$ is a torus, then $U/G$ can be taken to be a
product of projective spaces.  If $G = GL_n$, then $U/G$ can be taken
to be a Grassmannian (see the example in Section \ref{s.subset})

If $G$ is semisimple, then $U/G$ cannot be chosen
projective, for then the hyperplane class would be a nontorsion
element in $A^1_G$, but by Proposition \ref{conred},
$A^*_G \otimes \Q \cong S(\hat{T})^W \otimes
\Q$, which has no elements of degree 1.  Nevertheless for
semisimple (or reductive) groups we can obtain a cycle map
$$cl: A_*^G(X)_{\Q} \rightarrow H_{BM*}^T(X;\Q)^W$$ by identifying
$A_*^G(X) \otimes \Q$ with $A_*^T(X)^W \otimes \Q$ and
$H_{BM*}^G(X;\Q)$ with $H_{BM*}^T(X;\Q)^W$; if $X$ is compact then
the last group can be calculated with a compact model.

\section{Examples}  \label{examples}
In this section we calculate some examples of equivariant
Chow groups, particularly for connected groups. The point of this
is to show that computing {\it equivariant} Chow groups is no
more difficult than computing ordinary Chow groups, and in the case
of quotients, equivariant Chow theory gives a way of computing
ordinary Chow groups. Moreover,
since many of the varieties with computable Chow groups
(such as $G/P$'s, Schubert varieties,
spherical varieties, etc.) come with group actions, it is natural
to study their equivariant Chow groups.

\subsection{Representations and subsets} \label{s.subset}
For some groups there is a convenient choice of representations and
subsets.  In the simplest case, if $G= \G_m$ then we can take $V$ to an
$l$-dimensional representation with all weights $-1$, $U = V - \{0\}$, and
$U/G = \P^{l-1}$. If $G=T$ is a (split) torus of rank $n$, then we can
take $U = \oplus_1^n (V - \{0\})$ and $U/T = \Pi_1^n \P^{l-1}$. If $G =
GL_n$, take $V$ to be the vector space of $n \times p$ matrices ($p>n$),
with $GL_n$ acting by left multiplication, and let $U$ be the subset of
matrices of maximal rank.  Then $U/G$ is the Grassmannian $Gr(n,p)$.
Likewise, if $G=SL(n)$, then $U/G$ fibers over $Gr(n,p)$ as the
complement of the 0-section in the line bundle $\mbox{det}(S)
\rightarrow Gr(n,p)$, where $S$ is the tautological rank $n$ subbundle
on $Gr(n,p)$.

\subsection{Equivariant Chow rings of points} The equivariant Chow ring
of a point was introduced in \cite{To}.  If $G$ is connected
reductive, then $A^*_G \otimes \Q$ and (if $G$ is special) $A^*_G$ are
computed in \cite{E-G}.  The computation given there does not use a
particular choice of representations and subsets.  The result is that
$A^*_G \otimes \Q \cong S(\hat{T})^W \otimes \Q$, where $T$ is a maximal
torus of $G$, $S(\hat{T})$ the symmetric algebra on the group of
characters $\hat{T}$, and $W$ the Weyl group.  If $G$ is special this
result holds without tensoring with $\Q$.

\begin{prop} \label{conred}
Let $G$ be a connected reductive group with split maximal
torus $T$ and Weyl group $W$. Then
$A_*^G(X) \otimes \Q = A_*^T(X)^W \otimes \Q$. If $G$ is special
the isomorphism holds with integer coefficients.
\end{prop}

Proof. If $G$ acts freely on $U$, then so does
$T$. Thus for a sufficiently large representation $V$,
$A_{i}^T(X) = A_{i+l -t}((X \times U)/T)$ and
$A_i^G(X) = A_{i+l-g}((X \times U)/G)$ (here
$l$ is the dimension of $V$, $t$ the dimension of $T$
and $g$ the dimension of $G$). On the other hand,
$(X \times U)/T$ is a $G/T$ bundle over $(X \times U)/G$.
Thus
$A_{k}((X \times U/T)) \otimes \Q =A_{k+g-t}((X \times U)/G)^W \otimes \Q$
and if $G$ is special, then the equality holds integrally (\cite{E-G})
and the proposition follows.
\endproof \medskip

For $G$ equal to $\G_m$ or $GL(n)$ the choice of representations in
Section \ref{s.subset} makes it easy to compute $A^*_G$ directly,
without appealing to the result of \cite{E-G}.
If $l > i$, then $A^i_{\G_m} =  A^i(\P^{l-1}) = \Z \cdot t^i$, where
$t= c_1({\cal O}(1))$.  Thus, $A^*_{\G_m}(pt) = \Z[t]$.
More generally, for a torus of rank $n$, $A^*_{T}(pt) = \Z[t_1, \ldots , t_n]$.
Likewise, for $p$ sufficiently large, $A^*_{GL_{n}}(pt) =
A^i(Gr(n,p))$ is the free abelian group of homogeneous symmetric
polynomials of degree $i$ in $n$-variables (polynomials in the Chern
classes of the rank $n$ tautological subbundle).  Thus
$A^*_{GL_{n}}(pt) = \Z[c_1, \ldots , c_n]$ where $c_i$ has degree
$i$. Likewise $A^*_{SL_n}(pt) = \Z[c_2, \ldots c_n]$.

There is a map $A^*_{GL_{n}} \rightarrow A^*_{T}$, where $T$ is
a maximal torus.  This is a special case of a general construction:
if $G$ acts on $X$ and $H \subset G$ is a subgroup, then there is a
pullback $A_*^G(X) \rightarrow A_*^H(X)$. This map is induced by pulling
back along the flat map $X_H = X \times^H U \rightarrow X_G= X \times^G
U$.   We can identify the map $A^*_{GL_{n}} \rightarrow A^*_{T}$ concretely
as the map $\Z[c_1, \ldots c_n] \rightarrow \Z[t_1,
\ldots , t_n]$ given by $c_i \mapsto e_i(t_1, \ldots , t_n)$ (here
$e_i$ denotes the $i$-th symmetric polynomial), so $A^*_{GL_n}(pt)$ can
be identified with the subring of symmetric polynomials in $\Z[t_1,
\ldots , t_n]$.  This is a special case of the result of Proposition
\ref{conred}.

More elaborate computations are required to compute $A^*_G$ for other
reductive groups (if one does not tensor with $\Q$).  The cases $G =
O(n)$ and $G = SO(2n+1)$ have been worked out by Pandharipande
\cite{Pa2} and Totaro.  There is a conjectural answer for $G =
SO(2n)$, verified by Pandharipande for $n=2$.

\paragraph{Equivariant Chern classes over a point}  An equivariant
vector bundle over a point is a representation of $G$.  If $T = \G_m$,
equivariant line bundles correspond to the 1-dimensional
representation $L_a$ where $T$ acts by weight $a$.  If (as above) we
approximate $BT$ by $(V - \{0 \})/T = \P(V)$, where $T$ acts on $V$ with
all weights $-1$, then the tautological subbundle corresponds to the
representation $L_{-1}$.  Hence $c_T(L_a) = at$.

\subsection{Equivariant Chow rings of $\P^n$}
We calculate
$A^*_T(\P^n)$, where $T = \G_m$ acts diagonally on $\P^n$ with weights
$a_0, \ldots , a_n$ (i.e., $g \cdot (x_0:x_1 \ldots :x_n) = (g^{a_0}x_0:
g^{a_1}x_1 : \ldots : g^{a_n}x_n)$).  In this case, $X_T
\rightarrow U/T$ is the $\P^n$ bundle
$$
\P({\cal O}(a_0) \oplus \ldots \oplus {\cal O}(a_n)) \rightarrow \P^{l-1}.
$$
Thus $A^*(X_T) = A^*(\P^{l-1})[h]/(p(h,t))$ where $t$ is the
generator for $A^1(\P^{l-1})$ and
$$
p(h,t) = \sum_{i = 0}^{n} h^i e_i(a_0t, \ldots , a_n t).
$$   Letting the dimension of the
representation go to infinity we see that $A^*_T(\P^n) = \Z[t,h]/p(h,t)$.
Note that $A^*_T(\P^n)$ is a module of rank
$n+1$ over the $T$-equivariant Chow ring of a point.

Assume that the weights of the $T$-action are distinct.  Then the
fixed point set $(\P^n)^T$ consists of the points $p_0, \ldots , p_n,$
where $p_r \in \P^n$ is the point which is non-zero only in the $r$-th
coordinate. The inclusion $i_r: p_r \hookrightarrow \P^n$ is a regular
embedding.  The equivariant normal bundle is the equivariant vector
bundle over the point $p_r$ corresponding to the representation $V_r =
\oplus_{s \neq r} L_{a_s}$.  The equivariant pushforward $i_{r*}$ is
readily calculated.  For example, if $n = 1$ then $i_{r*}$ takes
$\alpha$ to $\alpha \cdot (h + a_s t)$ (where $s \neq r$).  Hence the
map $i_*: A_*^T((\P^1)^T) \rightarrow A^*_T(\P^1)$ becomes an
isomorphism after inverting $t$ (and tensoring with $\Q$ if $a_0$ and
$a_1$ are not relatively prime).  This is a special case of the
localization theorem for torus actions (\cite{EG38}).

We remark that the calculation of $A^*_T(\P^n)$ can be viewed as a special
case of the projective bundle theorem for equivariant Chow groups
(which follows from the projective bundle theorem for ordinary Chow
groups), since $\P^n$ is a projective bundle over a point, which is
trivial but not equivariantly trivial.

\subsection{Computing Chow rings of quotients}
By Theorem \ref{quotient} the rational Chow
groups of the quotient of a variety by a group acting with finite
stabilizers can be identified with the equivariant Chow groups of the original
variety.  If the original variety is smooth then the rational Chow
groups of the quotient inherit a canonical ring structure (Theorem
\ref{quotient.cor}).

For example, let $W$ be a representation of a split torus $T$ and let $X
\subset W$ be the open set on which $T$ acts properly. Since representations
of $T$ split into a direct sum of invariant lines,
it easy to show that $W - X$ is a finite union of invariant linear
subspaces $L_1,
\ldots ,L_r$.
When $T = {\bf G}_m$ and $X = W-\{0\}$ then the quotient is
a twisted projective space.

Let ${\cal R} \subset \hat{T}$ be the set of weights of $T$ on $V$.  If
$L \subset V$ is an invariant linear subspace, set $\chi_L = \Pi_{\chi
  \in {\cal R}} \chi^{d(L,{\chi})}$, where $d(L,{\chi})$ is the dimension
of the $\chi$-weight space of $V/L$.

\begin{prop} There is a ring isomorphism
$A^*(X/T)_{\Q} \simeq S(\hat{T})/(\chi_{L_1}, \ldots \chi_{L_r})$.
\end{prop}
Proof. By Theorem \ref{quotient}, $A^*(X/T)_{\Q} = A^T(X)_{\Q}$.  Since
$W-X$ is a union of linear subspaces $L_1, \ldots L_r$, we have an exact
sequence (ignoring the shifts in degrees)
$$
\oplus A^*_T(L_i) \stackrel{i_*} \rightarrow A^*_T(V) \rightarrow
A^*_T(U) \rightarrow 0 .
$$
Identifying $A^*_T(W)$ with $A^*_T(pt) = S(\hat{T})$, we see that
$A^T(U) = S(\hat{T})/im(i_*)$. Since each invariant linear
subspace is the intersection of invariant hypersurfaces,
the image of $A^*_T(L_i)$ in $A^*_T(W)=S(\hat{T})$
is the ideal $\chi^{c(L,\chi)}S(\hat{T})$.  \endproof \medskip

{\bf Remark}.
The preceding proposition is a simpler presentation of a computation in
\cite[Section 4]{Vi2}.  In addition, we do not need to assume that the
stabilizers are reduced, so there is no restriction on the characteristic.

In \cite{E-S1}, Ellingsrud and Str{\o}mme considered representations
$V$ of $G$ for which all $G$-semistable points are stable for a maximal
torus of $G$, and $G$ acts freely on the
set $V^s(G) $ of $G$-stable points. In this case they gave a
presentation for $A^*(V^s/G)$.
Using Theorem \ref{quotient}, it can be shown that their presentation
is valid (with $\Q$ coefficients) even if $G$ doesn't act freely on
$V^s(G)$.

In a more complicated example, Pandharipande (\cite{Pa1}) used
equivariant Chow groups to compute the rational Chow ring of the
moduli space, $M_{\P^r}(\P^k,d)$ of maps $\P^k \rightarrow \P^r$ of
degree $d$.  This moduli space is the quotient $U(k,r,d)/GL(k+1)$,
where $U(k,r,d) \subset \oplus^r_0 H^0(\P^k, {\cal O}_{\P^k}(d))$ is
the open set parameterizing base-point free $r+1$-tuples of polynomials
of degree $d$ on $\P^k$.  His result is that for any $d$,
$A^*(M_{\P^r}(\P^k,d))_{\Q}$ is canonically isomorphic to the rational
Chow ring of the Grassmannian $Gr(\P^k, \P^r)$.

\subsection{Intersecting equivariant cycles, an example} \label{noinprod}
Let $X = k^3 - \{ 0 \}$ and let $T$ denote the
$1$-dimensional torus acting with weights 1,2,2.  We let $A^i[X/T]$
denote the group of invariant cycles on $X$ of codimension $i$, modulo
the relation $\mbox{div}(f) = 0$, where $f$ is a $T$-invariant rational
function on an invariant subvariety of $X$.  We will show that there is
no (reasonable) intersection product, with integer coefficients, on
$A^*[X/T]$.  We will also compare $A^*[X/T]$ to $A^*_T(X)$, which does
have an integral intersection product.

Clearly $A^0[X/T] = \Z \cdot [X]$.  An invariant codimension $1$
subvariety is the zero set of a weighted
homogeneous polynomials $f(x,y,z)$, where $x$ has weight $1$ and $y$
and $z$ have weights $2$.  If $f$ has weight $n$ then the cycle defined
by $f$ is equivalent to the cycle $n \cdot p$, where $p$ is the class of
the plane $x = 0$.  Thus $A^1[X/T] = \Z \cdot p$.  The invariant
codimension $2$ subvarieties are just the $T$-orbits.  If we let $l$
denote the class of the line $x=y=0$, then we see that the orbit $T
\cdot (a,b,c)$ is equivalent to $l$ if $a=0$, and to $2l$ otherwise.
Thus $A^2[X/T] = \Z \cdot l$.  Finally, $A^i[X/T] = 0$ for $i \geq 3$.

If $Z_1$ and $Z_2$ are the cycles defined by $x=0$ and $y=0$, then $Z_1$
and $Z_2$ intersect transversely in the line $x=y=0$.  Thus, in a
``reasonable'' intersection product we would want $2 p^2 = [Z_1] \cdot
[Z_2] = l$ or $p \cdot p = \frac{1}{2}l$.  But $\frac{1}{2}l$ is not an
integral class in $A^2[X/T]$, so such an intersection product does not
exist.

Now consider the equivariant groups $A^*_T(X)$.  We model $BT$ by
$\P^N$, where $N$ is arbitrarily large; then the mixed space $X_T$
corresponds to the complement of the $0$-section in the vector bundle
${\cal O}(1) \oplus {\cal O}(2) \oplus {\cal O}(2)$.  Thus $A^*_T(X) =
\Z[t] / (4 t^3)$.

Each invariant cycle on $X$ defines an element of $A^*_T(X)$, so there
is a natural map $A^*[X/T] \rightarrow A^*_T(X)$.  This map takes $p$ to
$t$ and $l$ to $2 t^2$.  The equivariant theory includes the extra
cycle, $t^2$, necessary to define an integral intersection product.  We
can view elements of $A^*_T(X)$ as cycles on $X \times V$, where $V$ is
a representation of $T$ with all weights $1$.  The class $t^2$ is
represented by the cycle $x = 0, \phi = 0$ where $\phi$ is any linear
function on $V$.

\section{Intersection theory on quotients} \label{qint}
One of the most important properties of equivariant Chow groups is
that they compute the rational Chow groups of a quotient
by
a group acting with finite stabilizer. They can also be used to show
that the rational Chow groups of a moduli space which is a
quotient (by a group)
of a smooth algebraic space have an intersection product -- even
when there are infinitesimal automorphisms.

\subsection{Chow groups of quotients}
Let $G$ be a $g$-dimensional group acting on a algebraic space $X$.
Following Vistoli, we define a {\em quotient} $X \stackrel{\pi}
\rightarrow Y$ to be a map which satisfies the following properties
(cf. \cite[Definition 0.6(i - iii)]{GIT}): $\pi$ commutes with the
action of $G$, the geometric fibers of $\pi$ are the orbits of the
geometric points of $X$, and $\pi$ is submersive, i.e., $U \subset Y$
is open if and only if $\pi^{-1}(U)$ is.  (This is called a
topological quotient in \cite[Definition 2.7]{Kollar}.)  Unlike what
Mumford calls a
geometric quotient, we do not require that ${\cal O}_Y = \pi_*({\cal
O} _X)^G$. The advantage of this definition is that it is preserved under
base change.  In characteristic 0 there are no inseparable extensions,
so our quotient is in fact a geometric quotient (\cite[Prop. 0.2]{GIT}).
For proper actions, a geometric quotient is unique (\cite{GIT}, \cite{Kollar}).
If $L \subset K$ is an inseparable
extension and $G_K$ is a group defined
over $K$, then both $\mbox{Spec } L$ and $\mbox{Spec } K$
are quotients of $G_K$ by $G_K$. This example shows that
quotients need not be unique in characteristic $p$.


\begin{prop} \label{p.quotient}
(a) If $G$ acts with trivial stabilizer on an algebraic space $X$
and $X \rightarrow Y$ is the  principal bundle quotient then
$A_{i+g}^G(X) = A_{i}(Y)$.

(b) If in addition $X$ is quasi-projective and $G$ acts linearly,
then
$A_{i+g}^G(X,m) = A_i^G(X,m)$ for all $m \geq 0$.
\end{prop}
Proof. If the stabilizers are trivial,
then $(V \times X)/G$ is a vector bundle
over the quotient $Y$. Thus $X_G$ is an open set in this bundle
with arbitrarily high codimension, and the proposition follows from homotopy
properties of (higher) Chow groups. \endproof

\medskip

\begin{thm} \label{quotient}

(a) Let $X$ be an algebraic space with a (locally) proper $G$-action
and let
$X \stackrel{\pi} \rightarrow Y$ be a quotient.  Then
$$A_{i+g}^G(X) \otimes \Q \simeq A_i(Y) \otimes \Q.$$

(b) If in addition $X$ is quasi-projective with a linearized
$G$-action, and the quotient $Y$ is quasi-projective, then
$$A_{i+g}^G(X,m) \otimes \Q \simeq A_{i}^G(Y,m) \otimes \Q$$

\end{thm}

\begin{thm} \label{quotient.cor} With the same hypotheses
as in Theorem \ref{quotient}(a),
there is an isomorphism of operational Chow rings
$$\pi^*:A^*(Y)_{\Q} \stackrel{\simeq} \rightarrow A^*_G(X)_{\Q}.$$
Moreover if $X$ is smooth, then the map $A^*(Y)_{\Q} \stackrel{\cap
[Y]} \rightarrow A_*(Y)_{\Q}$ is an isomorphism.
In particular, if $X$ is smooth,
the rational Chow groups of the quotient space $Y=X/G$ have a ring
structure, which is independent of the presentation of $Y$ as a quotient
of $X$ by $G$.
\end{thm}

{\bf Remarks.}
(1) By \cite{Kollar} or \cite{KM}, if $G$ acts locally properly on a
(locally separated) algebraic space, then a (locally separated)
geometric quotient $X \stackrel{\pi} \rightarrow Y$ always exists in
the category of algebraic spaces.
Moreover, the result of \cite{KM} holds under the weaker
hypothesis that $G$ acts on $X$ with finite stabilizer.
However, the quotient $Y$ need not be locally separated.
We expect that Theorems \ref{quotient} and \ref{quotient.cor}
still hold in this case, but our proof does not go through.

(2) The hypotheses in Theorem \ref{quotient}(b) are purely
technical. They are necessary because the localization theorem for
higher Chow groups has only been proved for quasi-projective
schemes. If the localization theorem were proved for algebraic spaces,
Theorem \ref{quotient}(b) would hold in this case.

(3) Checking that an action is proper can be difficult. If $G$ is
reductive, then \cite[Proposition 0.8 and Converse 1.13]{GIT} give
criteria for properness when $G$ is reductive. In particular
if $X$ is contained in the set of stable points
for some linearized action of $G$ on $X$ then the action is proper.
If $X \rightarrow Y$ is
affine and $Y$ is quasi-projective the action is also proper. Not
surprisingly, checking that an action is locally proper is slightly
easier.  In particular if $G$ is a reductive and $X \rightarrow Y$ is
a geometric quotient such that $Y$ is a scheme, then the action is
locally proper if $X$ can be covered by invariant affine open sets.

(4) In practice, many interesting varieties arise as quotients of
smooth varieties by connected algebraic groups which act with finite
stabilizers.  Examples include simplicial toric varieties and various
moduli spaces such as curves, vector bundles, stable maps, etc.
Theorem \ref{quotient.cor} provides a tool to compute Chow groups of
such varieties (see Section \ref{examples} for some examples).

(5) As noted above, Theorem \ref{quotient.cor} shows that there
exists an intersection product on the rational Chow group of quotients
of smooth varieties and algebraic spaces.  There is a long history of
work on this problem.  In characteristic 0, Mumford \cite{Mu} proved
the existence of an intersection product on the rational Chow groups
of $\Mgbar$, the moduli space of stable curves.  Gillet \cite{Gi}
and Vistoli \cite {Vi} constructed intersection products on
quotients in arbitrary characteristic, provided that the stabilizers
of geometric points are reduced. (In characteristic 0 this condition
is automatic, but it can fail in positive characteristic.)  In
characteristic 0, Gillet (\cite[Thm 9.3]{Gi}) showed that his product
on $\Mgbar$ agreed with Mumford's, and in \cite[Lemma 1.1]{Ed} it was
shown that Vistoli's product also agreed with Mumford's.  If
the stabilizers are reduced, we show that our product agrees
with Gillet's and Vistoli's (Proposition \ref{triprod}).
Hence, Gillet's product and Vistoli's agree for quotient
stacks, answering a question in \cite{Vi}.
Moreover, Theorem
\ref{quotient.cor} does not require that the stabilizers be reduced
and is therefore true in arbitrary characteristic,
answering \cite[Conjecture 6.6]{Vi} affirmatively for
moduli spaces of quotient stacks.

(6) Equivariant intersection theory gives a nice way of working with
cycles on a singular moduli space ${\cal M}$ which is a quotient $X/G$
of a smooth variety by a group acting with finite stabilizer.  Given
a subvariety $W \subset {\cal M}$ and a family $Y
\stackrel{p}\rightarrow B$ of schemes parameterized by ${\cal M}$,
there is a map $B \stackrel{f} \rightarrow {\cal M}$. We wish to
define a class $f^*([W]) \in A_*B$ corresponding to how the image of
$B$ intersects $W$. This can be done (after tensoring with $\Q$) using
equivariant theory.

By Theorem \ref{quotient}, there is an isomorphism $A_*({\cal
M})_{\Q}= A_*^G(X)$ which takes $[W]$ to the equivariant class $w=
\frac{e_W}{i_W}[f^{-1}W]_G$.  Let $E \rightarrow B$ be the principal
$G$-bundle $B \times_{[X/G]} X$ (the fiber product is a scheme,
although the product is taken over the quotient stack $[X/G]$).
Typically, $E$ is the structure bundle of a projective bundle
$\P(p_*L)$ for a relatively very ample line bundle $L$ on $Y$).  Since
$X$ is smooth, there is an equivariant pullback $f^*_G: A_*^G(X)
\rightarrow A_*^G(E)$ of the induced map $E \stackrel{f_G}
\rightarrow X$, so we can define a class $f_G^*(w) \in
A_*^G(E)$. Identifying $A_*^G(E)$ with $A_*(B)$ we obtain our class
$f^*(W)$. When ${\cal M} = \Mgbar$ is the moduli space
of stable curves of genus $g$, these methods can be used
to re-derive formulas of \cite[Section 3]{Ed} for intersections
with various nodal loci.

\subsection{Preliminaries}
This section contains some results about quotients that
will be used in proving Theorem \ref{quotient}.  The reader may wish
to read the proofs after the proof of Theorem \ref{quotient}.

Let $G$ act locally properly on $X$ with quotient $X
\stackrel{\pi}\rightarrow Y$. The field extension $K(Y) \subset
K(X)^G$ is purely inseparable by \cite[p.43]{Borel}, and thus finite
because both $K(Y)$ and $K(X)^G$ are intermediate extensions of $k
\subset K(X)$ and $K(X)$ is a finitely generated extension of $k$. Set
$i_X = [K(X)^G:K(Y)]$.

Write $e_X$ for
the order of the stabilizer at a general point of $X$.  This
is the degree of the finite map $S(id_X) \rightarrow X$
where $S(id_X)$ is the stabilizer of the identity morphism
as defined in \cite[Definition 0.4]{GIT}. Note that the map
$S(id_X) \rightarrow X$ can be totally ramified. This occurs
exactly when the stabilizer of a general geometric point is non-reduced.
Finally,  set $\alpha_X = \frac{e_X}{i_X}$.

\begin{lemma} \label{l.factor}
Let $K = K(Y)$ be the ground field, and suppose $\pi: X \rightarrow Y
= \mbox{Spec }K$ is a quotient of a variety $X$ by a group
$G$ over $K$.  Then $\pi$ factors as $X \rightarrow \mbox{Spec }K(X)^G
\rightarrow Spec\; K(Y)$.
\end{lemma}
Proof. First, $X$ is normal.  To see this, let $Z \subset X$ be the
set of non-normal points, a proper $G$-invariant subset of $X$.  If
$L$ is an algebraically closed field containing $K$, write $X_L = X
\times_{Spec \; K} \mbox{Spec }L$, $G_L = G \times_{Spec \; K}
\mbox{Spec }L$.
Now, $Z_L$ is a proper $G$-invariant subset of
$X_L$.  Since $X_L$ is a single $G_L$-orbit, $Z_L$ is empty.  The map
$X_L \rightarrow X$ is surjective, by the going up theorem; hence $Z$
is empty, so $X$ is normal.

Now if $\mbox{Spec }A \subset X$ is an open affine subset, then $K(Y)
\subset A$,  $K(X)^G$ is integral over $K(Y)$, and $A$ is integrally
closed in $K(X)$.  We conclude that $K(X)^G \subset A$, which implies
the result.  \endproof

\medskip

{\bf Remark.}  The fact that $X_L$ is a single $G_L$-orbit is
essential to the result.  For example, suppose $K = {\bf F}_p(t)$, $A =
K[u,v]/(u^p - t v^p)$, $X = \mbox{Spec }A$, $G = \G_m$ acting by $g
\cdot(u,v) = (gu,gv)$.  The geometric points of $X_L$ form two
$G_L$-orbits, since $(X_L)_{red} = \A^1_L$.  The conclusion of the
lemma fails since $\frac{u}{v}$ is in $K(X)^G$ but not in $A$.

\bigskip

\begin{prop} \label{p.technical.quotient}
Let
$$
\begin{array}{ccc}
X' & \stackrel{g} \rightarrow &X \\
\small{\pi'} \downarrow & & \small{\pi} \downarrow\\
Y' & \stackrel{f} \rightarrow & Y
\end{array}
$$
be a commutative diagram of quotients with
$f$ and $g$ finite and surjective.  Then
$$
\frac{[K(X'):K(X)]}{[K(Y'):K(Y)]} = \frac{\alpha_X}{\alpha_{X'}}
= \frac{e_X}{e_{X'}} \cdot \frac{i_{X'}}{i_X}
$$
\end{prop}
Proof. Since we are checking degrees, we may replace $Y'$ and $Y$
by $K(Y')$ and $K(Y)$, and $X'$ and $X$ by their generic fibers over
$Y'$ and $Y$ respectively.
By the above lemma, we have a commutative diagram of varieties
$$\begin{array}{ccc}
X' &  \rightarrow & X \\
\downarrow & &  \downarrow\\
\mbox{spec}(K(X')^G) & \rightarrow & \mbox{spec}(K(X)^G)\\
\downarrow & &  \downarrow\\
\mbox{spec}(K(Y')) & \rightarrow & \mbox{spec}(K(Y)).
\end{array}$$
Since $i_{X'} := [K(X')^G:K(Y')]$ and $i_X :=[K(X)^G:K(Y)]$, it suffices
to prove that
$$
\frac{[K(X'):K(X)]}{[K(X')^G:K(X)^G]} = \frac{e_X}{e_{X'}}.
$$

By \cite[Prop. 2.4]{Borel} the extensions $K(X')^G \subset K(X')$ and
$K(X)^G \subset K(X)$ are separable (transcendental).  Thus, after
finite separable base extensions, we may assume that there are sections
$s':\mbox{spec}(K(X')^G) \rightarrow X'$ and $s: \mbox{spec}(K(X)^G)
\rightarrow X$.

The section $s$ gives us a finite surjective map $G_K \rightarrow X$,
where $G_K = G \times_{Spec \; k} \mbox{Spec }K(X)^G$.  The degree of this
map is $e_X = [K(G_K): K(X)]$ because $G_K \rightarrow X$
is the pullback of the action morphism $G \times X \rightarrow X$
via the map $X \rightarrow X \times X$ given by $x \mapsto (x,s(\mbox{Spec
}K))$.
Likewise, $e_{X'} = [K(G_{K'}): K(X)]$, where
$G_{K'} = G \times_{Spec \; k} \mbox{Spec }K(X')^G$.  Therefore,
$$
\frac{e_X}{e_{X'}} = \frac{[K(G_K): K(X)]}{[K(G_{K'}): K(X')]} =
 \frac{[K(X'):K(X)]}{[K(X')^G:K(X)^G]},
$$
since $[K(G_{K'}):K(G_K)] =
[K(X')^G:K(X)^G]$.  This completes the proof.  \endproof

\medskip

The following proposition is an analogue of \cite[Prop. 2.6]{Vi}
and \cite[Thm 6.1]{Seshadri}. Our proof is similar to Vistoli's.
\begin{prop} \label{whizzbang}
Suppose that $G$ acts locally properly on an algebraic space $X$.
Let $X
\rightarrow Y$ be a quotient. Then there is a commutative diagram of
quotients, with $X'$ a normal algebraic space:
$$\begin{array}{ccc} X' & \rightarrow  & X\\
\downarrow & & \downarrow\\
Y' & \rightarrow & Y
\end{array}$$
where  $X' \rightarrow Y'$ is a principal $G$-bundle (in particular
$G$ acts with trivial stabilizer on $X'$) and the horizontal
maps are finite and surjective.
\end{prop}

{\bf Remark.} If $X$ and $Y$ are quasi-projective, then so are $X'$
and $Y'$. If the action on $X$ is proper, then the action on
$X'$ is free.

Proof.
By \cite[Lemma, p.14]{GIT}, there is a finite map $Y' \rightarrow
Y$, with $Y'$ normal,
so that the pullback $X_1 \stackrel{\pi}\rightarrow Y'$ has a section
in a neighborhood of every point. Cover $Y'$
by a finite number of open sets $\{U_\alpha\}$ so that
$X_1 \rightarrow Y'$ has a section $U_\alpha \stackrel{s_{\alpha}}
\rightarrow
V_{\alpha}$ where $V_{\alpha} = \pi^{-1}(U_{\alpha})$.

Define a $G$-map
$$\phi_{\alpha}: G \times U_\alpha \rightarrow V_\alpha$$
by the Cartesian diagram
$$\begin{array}{ccc}
G \times U_\alpha & \stackrel{\phi_\alpha} \rightarrow & V_\alpha\\
\downarrow & & \small{id \times s_\alpha \circ \pi} \downarrow\\
G \times V_\alpha & \rightarrow & V_{\alpha} \times V_{\alpha}.
\end{array}$$
The action is locally proper so we can, by shrinking $V_\alpha$,
assume that
$\phi_\alpha$ is proper.
Since it is also quasi-finite, it is finite.

To construct a principal bundle $X' \rightarrow Y'$ we must glue
the $G \times U_{\alpha}$'s along
their fiber product over $X$. To do this we will find isomorphisms
$\phi_{\alpha\beta}: s_\alpha(U_{\alpha\beta}) \rightarrow
s_{\beta}(U_{\alpha\beta})$ which satisfy the cocycle
condition.

For each $\alpha, \beta$, let $I_{\alpha\beta}$ be the space
which parameterizes isomorphisms of $s_\alpha$ and $s_\beta$
over $U_{\alpha\beta}$ (i.e. a section $U_{\alpha\beta}
\rightarrow I_{\alpha\beta}$ corresponds to a global isomorphism
$s_\alpha(U_{\alpha\beta}) \rightarrow s_\beta(U_{\alpha\beta})$).
The space $I_{\alpha\beta}$ is finite
over $U_{\alpha\beta}$ (but possibly totally ramified
in characteristic
$p$) since it is defined by the cartesian diagram
$$\begin{array}{ccc}
I_{\alpha\beta} & \rightarrow & U_{\alpha\beta} \\
\downarrow & & \small{1 \times s_\beta} \downarrow\\
G\times U_{\alpha\beta} & \stackrel{1 \times \phi_\alpha} \rightarrow &
U_{\alpha\beta} \times V_{\alpha\beta}
\end{array}$$
(Note that $I_{\alpha\alpha}$ is the stabilizer of $s_\alpha(U_\alpha)$.)

Over $U_{\alpha\beta\gamma}$ there is a composition
$$I_{\alpha\beta} \times_{U_{\alpha\beta\gamma}} I_{\beta\gamma}
\rightarrow I_{\alpha\gamma}$$ which gives
multiplication morphisms which are surjective when $\gamma = \beta$.

After a suitable finite (but possibly inseparable) base change, we may
assume that there is a section $U_{\alpha\beta} \rightarrow
I_{\alpha\beta}$ for every irreducible component of $I_{\alpha\beta}$.
(Note that $I_{\alpha\beta}$ need not be reduced.) Fix an open set
$U_{\alpha}$. For $\beta \neq \alpha$ choose a section
$\phi_{\alpha\beta}: U_{\alpha\beta} \rightarrow I_{\alpha\beta}$.
Since the $I_{\alpha\beta}$'s split completely and $I_{\alpha\alpha}$
is a group over $U_{\alpha\alpha}$ (in the sense of \cite[Definition
0.1]{GIT}), there are sections
$\phi_{\beta\alpha}:U_{\alpha\beta} \rightarrow I_{\beta\alpha}$ so
that $\phi_{\alpha\beta} \cdot \phi_{\beta\alpha}$ is the identity
section of $U_{\alpha\alpha}$. For any $\beta, \gamma$ we can define a
section of $I_{\beta\gamma}$ over $U_{\alpha\beta\gamma}$ as the
composition $\phi_{\beta\alpha} \cdot \phi_{\alpha\gamma}$. Because
$I_{\beta\gamma}$ splits, the $\phi_{\beta\alpha}$'s extend to
sections over $U_{\beta\gamma}$.

By construction, the $\phi_{\beta\gamma}$'s satisfy the cocycle condition.
We can now define $X'$ by gluing the $G \times U_{\beta}$'s
along the $\phi_{\beta\gamma}$'s.  \endproof

\subsection{Proof of Theorems 3 and 4}
To simplify the notation we
give the proofs assuming that $G$ is connected (so the inverse image
in $X$ of a subvariety of $Y$ is irreducible).  All coefficients,
including those of cycle groups, are assumed to be rational.

\paragraph{Proof of Theorem \ref{quotient}}
Let $\Delta^m$ be the $m$-simplex of \cite{Bl}.
If $G$ acts locally properly on $X$, then
$G$ acts locally properly on $X \times \Delta^m$ by acting trivially on the
second factor. In this case, the boundary map of the higher Chow group
complex preserves invariant cycles, so there is a subcomplex of
invariant cycles $Z_*(X,\cdot)^G$.  Set
$$A_*([X/G],m) = H_m(Z_*(X, \cdot)^G,\partial).$$
This construction is well defined even when $X$ is an arbitrary
algebraic space.

Now if $X \rightarrow Y$ is a quotient, then so
is $X \times \Delta^m \stackrel{\pi} \rightarrow Y \times \Delta^m$.
Define a map
$\pi^*: Z_k(X,m) \otimes \Q \rightarrow Z_{k+g}(X,m)^G
\otimes \Q$ for all
$m$ as follows.  If $F \subset Y  \times \Delta^m$ is a $k+m$-dimensional
subvariety
intersecting the faces properly,
then $H = (\pi^{-1}F)_{red}$  is a $G$-invariant
$(k+m+g)$-dimensional subvariety of $X \times \Delta^m$ which intersects the
faces properly. Thus, $[H] \in Z_{k+g}(X,m)^G $.
Set $\pi^*[F] = \alpha_H[
H] \in Z_{k+g}^G(X,m)$, where $\alpha_H = \frac{e_H}{i_H}$ is defined
as above.
Since $G$-invariant subvarieties of $X \times \Delta^m$
correspond exactly to subvarieties of $Y \times \Delta^m$,
$\pi^*$ is an isomorphism of cycles for all $m$.

This pullback has good functorial properties:

\begin{prop} \label{p.functorial}
$(a)$ If
$$
\begin{array}{ccc}
X' & \stackrel{g} \rightarrow &X \\
\small{\pi'} \downarrow & & \small{\pi} \downarrow\\
Y' & \stackrel{f} \rightarrow & Y
\end{array}
$$
is a commutative diagram of quotients with $f$ and $g$ proper, then
$f_* {\pi '}^* = \pi^{*} g_*$ as maps $Z(Y',m) \rightarrow Z_*(X,m)^G$.

$(b)$ Suppose $T \subset X$ is a $G$-invariant subvariety. Let $S
\subset Y$ be its image under the quotient map. Set $U = X -T$ and
$V=U/G$, so there is a diagram of quotients:
$$
\begin{array}{ccc}
U & \stackrel{j} \rightarrow & X \\
\small{\pi} \downarrow & & \small{\pi} \downarrow \\
V & \stackrel{j} \rightarrow & Y .
\end{array}
$$
Then $\pi^*j^*  = j^* \pi^*$ as maps from $Z_k(Y,m)$ to $Z_{k+g}^G(U,m)$.

\end{prop}

Proof. Part (a) follows immediately from Proposition
\ref{p.technical.quotient}.  For (b), if $\alpha = [F]$ and $H=
\pi^{-1}(F)_{red}$, then $\pi^*j^*\alpha$ and $j^*\pi^*\alpha$ are
both multiples of $[H \cap U]$. Since $e_{[H \cap U]} = e_{[H]}$, and
$i_{[H \cap U]} = i_{[H]}$, the multiples are the same, proving
(b).  \endproof

\begin{prop} \label{squiggy}
$(a)$ The map $\pi^*$ commutes with the boundary operator defining higher
Chow groups.  In particular, there is an induced isomorphism of Chow
groups
$$A_k(Y,m) \simeq A_{k+g}([X/G],m)$$
(note that the higher Chow groups $A_k(Y,m)$ are defined as groups
even if $Y$ is only an algebraic space).

$(b)$ In the setting of Proposition \ref{p.functorial}(b), if
$X$ is quasi-projective with a linearized $G$-action, and
the quotient $Y$ is quasi-projective, then there is a
commutative diagram of higher Chow groups
$$
\begin{array}{ccccccc} \ldots \rightarrow & A_*([T/G],m) & \rightarrow &
A_*([X/G],m) & \rightarrow & A_*([U/G],m) & \rightarrow \ldots \\
& \uparrow & & \uparrow & & \uparrow & \\
\ldots \rightarrow & A_*(S,m) & \rightarrow &
A_*(Y,m) & \rightarrow & A_*(V,m) & \rightarrow \ldots
\end{array}$$
where the vertical maps are isomorphisms.  Hence the top row
of this diagram is exact.
\end{prop}

Proof. To prove (a), since $\pi^*$ is an isomorphism on the level of cycles,
once we
show that $\pi^*$ commutes with the boundary operator, it will follow
that the induced map on Chow groups is an isomorphism.  If
$$\begin{array}{ccc}
X' & \stackrel{g} \rightarrow & X \\
\small{\pi'} \downarrow & & \small{\pi} \downarrow\\
Y' & \stackrel{f} \rightarrow & Y
\end{array}$$
is a commutative diagram of quotients with $f$ and $g$ finite and
surjective, then $f_*$ and $g_*$ are surjective as maps of cycles. By
Proposition \ref{p.functorial}(a) it suffices to prove ${\pi '}^*:Z_*(Y')
\rightarrow Z_*(X')^G$ commutes with $\partial$.  By Proposition
\ref{whizzbang}, there exists such a
commutative diagram of quotients with $\pi':X'
\rightarrow Y'$ a principal bundle. Since $\pi'$ is flat,
${\pi '}^*$ commutes with $\partial$.  This proves (a).  Part (b)
follows from (a), Proposition \ref{p.functorial} and
the localization theorem for higher Chow groups.  \endproof

\medskip

Define a map
$\alpha:Z_*([X/G],m) \rightarrow Z_*(X_G,m)$
by $[F]  \mapsto [F]_G$.
This map commutes with proper pushforward and flat pullback.
Arguing as in Proposition \ref{squiggy}(a), we see that $\alpha$
commutes with the boundary operator defining higher Chow groups
and hence induces a map on Chow groups (again denoted $\alpha$).

The
proof of Theorem \ref{quotient} is completed by the following
proposition.

\begin{prop} \label{warhol}
(a) If $X$ is an algebraic space with a (locally) proper $G$-action, then the
map $\alpha: A_*([X/G]) \rightarrow A^G_*(X)$ is an isomorphism.

(b) If $X$ is quasi-projective with a proper linearized $G$-action and
a quasi-projective quotient $X \rightarrow Y$ exists, then
$\alpha: A_*([X/G],m) \rightarrow A^G_*(X,m)$ is an isomorphism
for $m > 0$.
\end{prop}

Proof of (a).
Let $Y$ be the quotient of $X$ by $G$.
We will prove that the composition
$\alpha \circ \pi^*:A_*(Y) \rightarrow A_*^G(X)$ is an isomorphism.
By Proposition \ref{whizzbang} there is
a commutative diagram of quotients with $g$ and $f$ finite
and surjective
$$
\begin{array}{ccc}
X' & \stackrel{g} \rightarrow &X \\
\small{\pi'} \downarrow & & \small{\pi} \downarrow\\
Y' & \stackrel{f} \rightarrow & Y
\end{array}
$$ and $X' \rightarrow Y'$ is a principal bundle.

Since $X' \rightarrow X$ and $Y' \rightarrow Y$
are finite and surjective,
\cite[Theorem 1.8]{Kimura} (which extends to the equivariant setting)
says that there are exact sequences
$$\begin{array}{c}
A_*^G(X' \times_X X') \stackrel{g_{1*} -
g_{2*}}\rightarrow A_*^G(X') \stackrel{g_*} \rightarrow
A_*^G(X) \rightarrow 0\\
A_*^G(Y' \times_Y Y') \stackrel{f_{1*} -
f_{2*}}\rightarrow A_*^G(Y') \stackrel{f_*} \rightarrow
A_*^G(Y) \rightarrow 0
\end{array}$$
where $g_i$ and $f_i$ are the projections to $X'$ and $Y'$.
Set $X'' = X' \times_X X'$, and set $Y'' = X''/G$.
The natural map $p:Y'' \rightarrow (Y' \times_Y Y')$
is finite and surjective, so the pushforward $A_*(Y'') \rightarrow
A_*(Y' \times_Y Y')$ is a surjection.  Hence the second sequence
remains exact if we replace $Y' \times_Y Y'$ by $Y''$ (and $f_i$ by
$f_i \circ p$).

We have a commutative diagram of exact sequences
$$\begin{array}{cccccc}
A_*^G(X'') &  \rightarrow & A_*^G(X') & \rightarrow &
A_*^G(X) & \rightarrow 0\\
\uparrow & &  \uparrow & & \uparrow  &  \\
A_*(Y'') & \rightarrow & A_*(Y') & \rightarrow &
A_*(Y) & \rightarrow 0 ,
\end{array}$$
where the vertical maps are $\pi^{''*} \circ \alpha{''}$, $\pi^{'*}
\circ \alpha{'}$, and $\pi \circ \alpha$, respectively.  By
Proposition \ref{p.quotient}, the first two maps are isomorphisms, so
by the 5-lemma, the third is as well.

\medskip

Proof of (b): We are going to use the localization
exact sequences for higher equivariant Chow groups, and for the invariant
Chow groups $A_*([X/G],m)$ (Proposition \ref{squiggy}(b)).  This is why
we must assume that $X$ is quasi-projective and that a
quasi-projective quotient exists.

By Proposition \ref{whizzbang}, we can find a finite surjective map
$g: X' \rightarrow X$, where $G$ acts freely on $X'$.  The map
$\alpha': A_*([X'/G],m) \rightarrow A^G_*(X',m)$ is an isomorphism, as
noted in (a).  By Noetherian induction and the localization long exact
sequences it suffices to prove the result when $X$ is replaced by the
open subset over which $g$ is flat, so assume this.  Because $g$ is
also finite,  $g_* g^*$ is multiplication by the degree of
$g$. Hence (since we are using rational
coefficients) the flat pullback $g^*$ makes $A_*^G(X,m)$ a summand in
$A_*^G(X',m)$ and $A_*([X/G],m)$ a summand in $A_*([X'/G],m)$.  Since
$\alpha' \circ g^* = g^* \circ \alpha$, the summand $A_*^G([X/G],m) \subset
A_*^G( [X'/G],m)$ is isomorphic to the summand $A_*^G(X,m) \subset
A_*^G([X/G],m)$.  This proves the proposition, and with it Theorem
\ref{quotient}.  \endproof.

\medskip

\paragraph{Proof of Theorem \ref{quotient.cor}}  The proof is similar to
\cite[Proposition 6.1]{Vi}. Let
$\pi: X \rightarrow Y$ be the quotient map.
We define a pullback $\pi^*:A^*(Y) \rightarrow A^*_G(X)$
as follows:
Suppose $c \in A^i(Y)$, $Z
\rightarrow X$ is a $G$-equivariant morphism, and $z \in
A_*^G(Z)$.
For any representation, there are maps $Z_G \rightarrow X_G \rightarrow Y$.
The class $z$ is represented by a class $z_V \in A_{*+l-g}(Z_G)$
for some mixed space $Z_G$. Define
$$\pi^*c \cap z = c \cap z_V \in A_{*+l-g-i}(Z_G) \simeq A_{*-i}^G(Z) $$
As usual, this definition is independent of the representation,
so $\pi^*c \cap \alpha \in A_*^G(Z)$.

Let $\hat{\pi}: A_* (Y) \rightarrow A_{*+g}^G(X)$ denote the
isomorphism of Theorem \ref{quotient}.

\begin{lemma} \label{l.compatibility}
If $c \in A^*(Y)$ and $y \in A_*(Y)$, then
$$\hat{\pi} (c \cap y) = \pi^*c \cap \hat{\pi} y.$$
\end{lemma}
Proof of Lemma \ref{l.compatibility}.
We first prove the lemma when $X \rightarrow Y$ is a principal
$G$-bundle. In this case the map $\hat{\pi}: A_*(Y) \rightarrow
A_{*+g}^G(X) \simeq A_{* + g}(X_G)$ is just the
pullback induced by flat map $X_G \rightarrow Y$.
Since the operations in $A^*(Y)$ are compatible with flat pullback
the lemma follows in this case.

By proposition \ref{whizzbang} there is a commutative diagram
of quotients with $g$ and $f$ finite and surjective.
$$
\begin{array}{ccc}
X' & \stackrel{g} \rightarrow &X \\
\small{\pi'} \downarrow & & \small{\pi} \downarrow\\
Y' & \stackrel{f} \rightarrow & Y
\end{array}
$$ and $X' \rightarrow Y'$ is a principal bundle.
Then $y = f_*(y')$ for some $y' \in A_*(Y')$.
Since $c  \cap y' = f^*c \cap y'$, we have by the first case
$$\hat{\pi'}(c \cap y') = \pi'^*f^*c \cap \hat{\pi'}y'$$
Let $x'_V \in A_*^G(X'_GXS)$ be the class corresponding to $\hat{\pi'}y'$.
Then $\pi'^*f^*c \cap \hat{\pi'}y' = c \cap x'_V = g^*\pi^*c \cap
\hat{\pi'}y$.
Thus we obtain the equation
$$(*) \mbox{  }\hat{\pi'}(c \cap y') =g^*\pi^*c \cap \hat{\pi'}y.$$
By Proposition \ref{p.functorial} $g_*(\hat{\pi'}(c \cap y'))
= \hat{\pi}(c \cap f_*y') = \hat{\pi}(c \cap y)$. The lemma
follows by applying $g_*$ to both sides of (*). \endproof

Given the lemma, we show that $\pi^*$ is an isomorphism as follows.
For injectivity, it suffices (by base change) to show that if $\pi^* c
= 0$ then $c \cap y = 0$ for all $y \in A_*(Y)$; this follows since
$\hat{\pi} (c \cap y) = 0$ by the lemma, and $\hat{\pi}$ is an
isomorphism.

The proof of surjectivity is more subtle.

Given  $d \in A^*_G(X)$ define
$c \in A^*(Y)$ as follows:  If $Y' \rightarrow Y$ and $y' in
A_*(Y')$, set
$$
c \cap y' =  \hat{\pi}'^{-1} (d \cap \hat{\pi}' y).
$$ where $\pi': X' = X \times_Y Y' \rightarrow Y'$ is
the pullback quotient.

We must now show $\pi^*c = d$.

We begin with a preliminary construction.
For any mixed space $X_G$
define a map $r:A^i_G(X) \rightarrow A^i(X_G)$
as follows: Given $Z \rightarrow X_G$ let $Z_U \rightarrow Z$
be the pullback of the principal $G$-bundle $X \times U
\rightarrow X_G$. Since $Z_U \rightarrow  Z$ is
a principal bundle we can identify  $A_*(Z)$ with $A_{* +g}^G(Z_U)$
and view a class $z \in A_*(Z)$ as an equivariant class in $A_{*+g}^G(Z_U)$.
Now set $r(c) \cap z = c \cap z \in A_*(Z_U) \simeq A_*(Z)$.
{}From the construction it is clear that if $X' \rightarrow X$ is equivariant
and $x' \in A_*^G(X')$ corresponds to $x'_V \in A_*(X'_G)$
then $r(c) \cap x'_V$ corresponds to  $c \cap x'$
under the identification of $A_*^G(X')$ $A_*^G(X'_G)$.

Suppose $Z \rightarrow X$ is equivariant
and $z \in A_*^G(Z)$ is defined by the class $z_V \in A_*^G(Z_G)$.
Let $\pi':X' \rightarrow Z_G$ be the quotient induced by pulling
back $X \rightarrow Y$  along the map $Z_G \rightarrow Y$.
By definition we must show $\hat{\pi}'^{-1}(d \cap \hat{\pi}'z_V)
\in A_*(Z_G) \simeq A_*^G(Z)$ defines the same class as $d \cap z
\in A_*^G(Z)$. By construction, $d \cap z$ is defined by
the class $r(d) \cap z_V \in A_*(Z_G)$. Since
$\hat{\pi}'$ is an isomorphism it suffices to prove
that $d \cap \hat{\pi}'z_V$ equals $\hat{\pi}'(r(d) \cap z_V)$.
By Lemma \ref{l.compatibility}
$\hat{\pi}'(r(d) \cap z_V) = \pi'^{*}r(d) \cap \hat{\pi}'z_V$.
The class $\hat{\pi}'z_V$ is represented by a class $(\hat{\pi}'z_V)_V
\in A_*(X'_G)$ and
$\pi^{'*}r(d) \cap \hat{\pi}'z_V = r(d) \cap (\hat{\pi}'z_V)_V$.
Finally $r(d) \cap (\hat{\pi}'z_V)_V = d \cap \hat{\pi}'z_V$
under the identification of $A_*(X'_G)$ with $A_*^G(X')$.
This proves that $\pi^*$ is surjective and with it part (a)
of the theorem.

To prove (b)
recall that  $\hat{\pi}(c \cap y) = \pi^*c \cap \hat{\pi}y$.
By Proposition \ref{opsmooth} the map $A^*_G(X) \stackrel{\cap [X]_G}
\rightarrow A_*(X)$ is an isomorphism (with $\Z$ coefficients).
Since $\hat{\pi}([Y]) = \alpha_X [X]_G$, the map
the map $c \mapsto c \cap [Y]$ is an isomorphism (with $\Q$ coefficients).
\endproof.

\section{Intersection theory on quotient stacks and their moduli}
\label{itstacks}
If $G$ acts on an algebraic space $X$, a quotient $[X/G]$ exists in
the category of stacks (\cite[Example 4.8]{D-M}; see below).
This section relates equivariant
Chow groups to Chow groups of quotient stacks.

We show that for proper actions, with rational coefficients,
equivariant Chow groups coincide with the Chow groups defined
by Gillet in terms of integral substacks.  Thus, in this case
the intersection products of Gillet and Vistoli are the same.
For an arbitrary action, the equivariant Chow groups are
an invariant of the quotient stack.  As an application, we calculate
the Chow rings of the moduli stacks of elliptic curves.

\subsection {Definition of quotient stacks}
We recall the definition of quotient stack. For an introduction
to stacks see \cite{D-M}, \cite{Vi}.
Let $G$ be a linear algebraic group acting
on a scheme (or algebraic space). The quotient stack
$[X/G]$ is the stack associated to
the groupoid $G \times X \rightarrow X \times X$.
If $B$ is a scheme, sections
of $[X/G](B)$ are principal $G$-bundles $E \rightarrow B$
together with an equivariant map $E \rightarrow X$.
Morphisms in $[X/G](B)$ are $G$-bundle isomorphisms
which preserve the map to $X$. In particular, if
$G$ acts on $X$ with trivial stabilizers then
all morphisms are trivial, $[X/G]$ is
a sheaf in the \'etale topology and the
quotient is in fact an algebraic space (Proposition \ref{l.algspacequotient}).

If $G$ acts with finite reduced stabilizers then the diagonal
$[X/G] \rightarrow [X/G] \times [X/G]$ is unramified
and $[X/G]$ is a Deligne-Mumford stack. If
$G$ acts (locally) properly then the diagonal
$[X/G] \rightarrow [X/G] \times [X/G]$ is (locally)
proper and the stack is (locally) separated. In
characteristic $0$ any (locally) separated stack
is Deligne-Mumford. However, in characteristic $p$
this need not be true since the stabilizers of geometric
points can be non-reduced.

\subsection{Quotient stacks for proper actions}
If $G$ acts locally properly on $X$ with reduced stabilizers then
the quotient stack $[X/G]$ is locally
separated and Deligne-Mumford. The rational Chow groups
$A_*([X/G]) \otimes \Q$ were first defined by Gillet \cite{Gi} and
coincide with the groups $A_*([X/G]) \otimes \Q$ defined above.  More
generally, if $G$ acts with finite stabilizers which are not reduced,
then Gillet's definition can be extended and we can define the
``naive'' Chow groups $A_k([X/G])_{\Q}$ as the group generated by
$k$-dimensional integral substacks modulo rational equivalences.
In this context,
Proposition \ref{warhol} can be restated in the
language of stacks as
\begin{prop}
Let $G$ be a $g$-dimensional group which acts
locally properly on an algebraic space $X$
(so the quotient $[X/G]$ is a locally separated Artin stack).
Then $A_i^G(X) \otimes \Q = A_{i-g}([X/G]) \otimes \Q$.
\endproof
\end{prop}

{\bf Remark.} Although $A_*^G(X) \otimes \Q = A_*([X/G]) \otimes \Q$,
the integral Chow groups may have non-zero torsion for all $i < \mbox{dim }
X$.

\medskip

With the identification of $A_*^G(X) \otimes \Q$ and $A_*([X/G]) \otimes \Q$
there are three intersection products on the rational Chow groups
of a smooth Deligne-Mumford quotient stack
-- the equivariant product, Vistoli's product
defined using a Gysin pullback for regular embeddings of stacks,
and Gillet's product defined using the product in higher $K$-theory.
The next proposition shows that they are identical.

\begin{prop} \label{triprod}
If $X$ is smooth and $[X/G]$ is a locally separated
Deligne-Mumford stack
(so $G$ acts with finite, reduced stabilizers)
then
the intersection products on $A_*([X/G])_{\Q}$ defined by Vistoli and
Gillet are the same as the equivariant product on $A_*^G(X)_{\Q}$.
\end{prop}

Proof. If $V$ is an $l$-dimensional representation, then all three
products agree on the smooth quotient space $(X \times U)/G$
(\cite{Vi}, \cite{Grayson}). Since the flat pullback of stacks
$f:A^*([X/G])_{\Q} \rightarrow A^*((X \times U)/G)_{\Q}$ commutes with
all 3 products, and is an isomorphism up to arbitrarily high codimension,
the proposition follows.  \endproof

\subsection{Integral Chow groups of quotient stacks} \label{intstack}
Suppose $G$ acts arbitrarily on $X$. Consider the (possibly non-separated
Artin) quotient
stack $[X/G]$.
The next proposition shows that the equivariant Chow groups
do not depend on the presentation as a quotient, so they are an invariant
of the stack.

\begin{prop} \label{qstacks}
Suppose that $[X/G] \simeq [Y/H]$ as quotient stacks. Then
$A_{i+g}^G(X) \simeq A_{i+h}^H(Y)$, where $\mbox{dim } G = g$
and $\mbox{dim } H = h$.
\end{prop}
Proof. Let
$V_1$ be an $l$-dimensional representation of $G$, and $V_2$ an
$M$ dimensional representation of $H$. Let $X_G = X \times^G U_1$
and $Y_H = X \times^H U_2$, where $U_1$ (resp. $U_2$) is an open set
on which $G$ (resp. $H$) acts freely.
Since the diagonal of a quotient stack is representable,
the fiber product $Z=X_G \times_{[X/G]} Y_H$ is an algebraic space. This
space is a bundle over $X_G$ and $Y_H$ with fiber $U_2$ and $U_1$
respectively.
Thus, $A_{i+l}(X_G) = A_{i+l+m}(Z) = A_{i+m}(Y_H)$
and the proposition follows.
\endproof

\medskip

As a consequence of Proposition \ref{qstacks} we can define the
integral Chow groups of a quotient stack ${\cal F} = [X/G]$ by
$A_i({\cal F}) = A_{i-g}^G(X)$ where $g = \mbox{dim }G$.

\begin{prop} If ${\cal F}$ is smooth, then $\oplus A_*({\cal F})$ has an
integral
ring structure.
\endproof \end{prop}

Following \cite[Definition, p. 64]{Mu} we define the
Picard group $Pic_{fun}({\cal F})$ of an algebraic stack ${\cal
F}$ as follows.  A element ${\cal L} \in Pic_{fun}( {\cal
F})$ assigns to any map $S \stackrel{F} \rightarrow {\cal F}$ of a
scheme $S$, an isomorphism class of a line bundle $L(F)$ on
$S$. Moreover the assignment must satisfy the following compatibility
conditions.

(i)  Let $S_1 \stackrel{F_1} \rightarrow {\cal F}$, $S_2
\stackrel{F_2} \rightarrow {\cal F}$ and $S_2 \stackrel{F_2} \rightarrow
{\cal F}$ be maps of schemes to ${\cal F}$.
If there is a map $S_1 \stackrel{f}
\rightarrow S_2$ such that $F_1 = F_2 \circ f$ then there
is an isomorphism $\phi(f): L(F_1) \simeq f^*(L(F_2))$.

(ii) If $S_1 \stackrel{f} \rightarrow S_2 \stackrel{g} \rightarrow S_3$
are maps of schemes such that $F_2 = F_3 \circ g$ and $F_1 = F_2 \circ f$
then there is a commutative diagram of isomorphisms
$$\begin{array}{ccc}
L(F_1) & \stackrel{\phi(F_1)} \rightarrow & f^*(L(F_2)) \\
\small{\phi(g \circ f)} \downarrow & & \small{f^*\phi(g)} \downarrow\\
(g \circ f)^*(L(F_3)) & = & f^*(g^*(L(F_3))).
\end{array}$$
The product ${\cal L} \otimes {\cal M}$ assigns
to $S \stackrel{f} \rightarrow {\cal F}$ the line
bundle $L_f \otimes M_f$.

\begin{prop}
Let $X$ be a smooth variety with a $G$-action. Then $A^1_G(X) =
Pic_{fun}([X/G])$.
\end{prop}
Proof. By Theorem \ref{piciscool}, $A^1_G(X) = Pic^G(X)$. The latter group
is naturally isomorphic to $Pic_{fun}([X/G])$.
\endproof \medskip

More generally, if ${\cal F}$ is any stack, we can define the integral
operational Chow ring $A^*({\cal F})$ as follows.  An element $c \in
A^k({\cal F})$ defines an operation $A_*B \stackrel{c_f} \rightarrow
A_{*-k}B$ for any map of a scheme $B \stackrel{f} \rightarrow {\cal
F}$.  The operations should be compatible with proper
pushforward, flat pullback and intersection products for maps of
schemes to ${\cal F}$ (cf. \cite[Definition 17.1]{Fulton} and
\cite[Definition 5.1]{Vi}). (This definition differs slightly from
Vistoli's because we use integer coefficients and
only consider compatibility with maps of schemes.)

\begin{prop}
Let ${\cal F} \simeq [X/G]$ be a smooth quotient stack. Then
$A^*({\cal F}) = A^*_G(X)$.
\end{prop}
Proof. Giving a map $B \stackrel{f} \rightarrow [X/G]$ is equivalent
to giving a principal $G$-bundle $E \rightarrow B$ together with
an equivariant map $E \rightarrow X$. An element of $A_G^*(X)$ defines
an operation on $A_*^G(E) = A_*(B)$, hence an operational
class in $A^*({\cal F})$. Conversely an operational class
$c \in A^k({\cal F})$ defines an operation on $A_*(X_G)$
corresponding to the map $X_G \rightarrow {\cal F}$ associated
to the principal bundle $X \times U \rightarrow X_G$. Set
$d = c \cap [X_G] \in A_{dim \; X -k}^G(X)$.
Since $X$ is smooth, the latter group is isomorphic to $A^k_G(X)$.
\endproof

\medskip

{\bf Remark.} Proposition \ref{qstacks} suggests that there should be
a notion of Chow groups of an arbitrary algebraic stack which can be
non-zero in arbitrarily high degree. This situation would be
analogous to the cohomology of quasi-coherent sheaves on the \'etale
(or flat) site (cf. \cite[p. 101]{D-M}).

\subsection{The Chow ring of the moduli stack of elliptic curves}
In this section
we will work over a field of
characteristic not equal to 2 or 3, and compute the Chow ring of the moduli
stacks ${\cal M}_{1,1}$ and $\overline{{\cal M}}_{1,1}$ of elliptic curves.
A. Vistoli has independently
obtained the results of this section, also using equivariant intersection
theory.  He has also calculated the Chow ring of ${\cal M}_2$, the moduli
stack of elliptic curves. This calculation will appear as an appendix
to this article \cite{Viap}.

\paragraph{Construction of the moduli stack}
The stacks ${\cal M}_{1,1}$ and $\overline{{\cal M}}_{1,1}$ are
defined as follows.  A section of ${\cal M}_{1,1}$
(resp. $\overline{{\cal M}}_{1,1}$) over $S$ is a family $(X
\stackrel{\pi} \rightarrow S, \sigma)$ where $X \rightarrow S$ is a
smooth (resp. possibly nodal) curve of genus $1$ and $\sigma: S \rightarrow
X$ is a smooth section.

Our construction of ${\cal M}_{1,1}$ and $\overline{{\cal M}}_{1,1}$
is similar to \cite[Section 4]{Mu}.
Let $\P(V) = \P^9$ be the projective space
of homogeneous degree $3$ forms in variables $x, y, z$. Let
$X \simeq \A^3 \subset \P(V)$
be the affine subspace parameterizing forms proportional to
$$y^2 z - (x^3 + e_1 x^2 z + e_2 x z^2 + e_3 z^3),$$ with $e_1, e_2,
e_3$ arbitrary elements of the field $k$.  Let $G = \{ \left(
\begin{array}{ccc} a & 0 & b \\ 0 & c & 0 \\ 0 & 0 & d \end{array}
\right) | a^3 = c^2d \neq 0\}$. The image of $G$ in  $PGL(3)$
consists of projective transformations which stabilize $X$.

Since
$a \neq 0$ for all $g \in G$, we can normalize and identify $G$ with
the subgroup\\
$\{ \left( \begin{array}{ccc} 1 & 0 & B \\ 0 & A & 0 \\ 0
& 0 & A^{-2} \end{array} \right) | A \neq 0 \}$ of $GL(3)$.  An
element $g = \left( \begin{array}{ccc} 1& 0 & B \\ 0 & A & 0 \\ 0 & 0
& A^{-2} \end{array} \right)$ acts on $(e_1, e_2, e_3)$ by
$$(e_1, e_2, e_3) \mapsto (A^{-2} e_1 + 3B, A^{-4} e_2 + 2 A^{-2} B e_1 +
3B^2, A^{-6} e_3 + A^{-4} B e_2 + A^{-2} B^2 e_1 + B^3).$$

Let $U \subset X$ be the open set where the polynomial
$(x^3 + e_1 x^2  + e_2 x
+ e_3)$ has distinct roots over the algebraic closure of the field
of definition. Likewise, let $W \subset X$ be the open set
where $x^3 + e_1 x^2 + e_2 x$ has at least 2 distinct roots.

\begin{prop}
(a) ${\cal M}_{1,1} \simeq  [U/G]$.

(b) $\Mbar_{1,1} \simeq [W/G]$.
\end{prop}
Proof of (a). Let $Y^{\cdot}$ be the functor such that a section of
$Y^{\cdot}$ over $S'$ is a
triple $(X' \stackrel{f'} \rightarrow S', \sigma', \phi')$ where $(X',\sigma')$
is an elliptic curve with section $\sigma'$ and $\phi'$
is an isomorphism of the flag ${\cal O}_{S'} \subset {\cal O}_{S'}^{\oplus 2}
\subset {\cal O}_{S'}^{\oplus 3}$ with the flag $f'_*({\cal O}_{X'}(\sigma))
\subset f'_*({\cal O}_{X'}(\sigma)^{\otimes 2}) \subset
f'_*({\cal O}_{X'}(\sigma)^{\otimes 3})$.  This functor is represented by the
scheme $Y$ parameterizing cubic forms in $3$ variables with a flex
at $(1:0:0)$, since there is a universal triple
$({\cal X} \stackrel{F} \rightarrow Y, \sigma_Y, \Phi)$ on $Y$.
There is an obvious
action of the group $B \subset GL(3)$ of upper triangular
matrices on $Y$.  Claim: The quotient stack $[Y/B]$ is ${\cal M}_{1,1}$.
Proof: Given a family $X \stackrel{f} \rightarrow S$ of elliptic curves,
there is a canonical principal $B$-bundle $S' \rightarrow S$
associated to the flag of vector bundles
${\cal O}_S \subset {\cal O}_S^{\oplus 2}
\subset {\cal O}_S^{\oplus 3}$.  By construction the pullback of
this flag to $S'$ is equipped with an isomorphism with the standard flag
in ${\cal O}_{S'}^{\oplus 3}$.  Thus from a family
of elliptic curves over $S$, we obtain a principal
$B$-bundle $S' \rightarrow S$ together with an
equivariant map $S' \rightarrow Y$, i.e., a section of $[Y/B]$ over
$S$.  The converse is similar.

Thus, it suffices to prove that $[U/G] \simeq [Y/B]$.
Let $Y'$ be the quotient of $Y$ by the natural $\G_m$-action.
Then $[Y/B] = [Y'/B']$, where $B'$ is the group of upper-triangular
matrices in $PGL(3)$. Identify $Y'$ with the locally
closed subscheme of $\P^9$ corresponding to cubics
with a flex at $(1:0:0)$.
Identifying $G \subset B'$, the inclusion $U \subset Y'$
is $G$-equivariant for the corresponding actions of $G$ and $B'$.
Thus we obtain a smooth (representable) morphisms of
quotient stacks
$$[U/G] \rightarrow [Y'/G] \rightarrow [Y'/B']
\simeq  {\cal M}_{1,1}.$$
Two curves in $U$ are isomorphic iff they are in the same $G$-orbit,
so the map $[U/G] \rightarrow {\cal M}_{1,1}$
is quasi-finite. Moreover, every elliptic curve
can be embedded in $\P^2$ as double cover of $\P^1$ branched
at $\infty$ and 3 finite points. Thus the map is surjective.
Finally, a direct check shows that the stabilizer of
the $G$-action on a point of $U$ is exactly the automorphism
group of the corresponding elliptic curve. Thus the geometric
fibers are single points. Hence, the map is \'etale, surjective
of degree 1, and thus an isomorphism of
stacks.

The proof of (b) is essentially identical.
\endproof

\medskip

We now compute the Chow ring of this stack.

\begin{prop} \label{didwediss}
(a) $A^*({\cal M}_{1,1}) = \Z[t]/12t$.

(b) $A^*(\overline{{\cal M}}_{1,1})= \Z[t]/24t^2$.
\end{prop}
Proof.
(a) By definition,  $A^*({\cal M}_{1,1}) = A^*_G(U)$.
The equivariant Chow ring
is not hard to calculate. Let $T \subset G$ be the
maximal torus. Then $G$ is a unipotent extension of $T$, so
$A^*_T(U) = A^*_G(U)$. Now $T$ acts diagonally on $X = \{(e_1,e_2,e_3)\}$
with weights $(-2,-4,-6)$. Let $S = X - U$, then $S$ is the discriminant
locus in $\A^3$, which
can be identified as the image of the big diagonal (i.e., the image
under the $S_3$ quotient map $A^3 \rightarrow X \simeq A^3$
which maps $(a,b,c) \mapsto (a+b+c, ab + bc + ac, abc))$ and has
equation $f(e_1, e_2, e_3) = 4e_2^3 + 27 e_3^2 - 18e_1 e_2 e_3 - e_1^2
e_2^2 + 4e_1^3e_3$.
The form $f$ is homogeneous of weighted degree $-12$ with respect to
the $T$-action on $X$.
Since $S \subset X$ is a divisor, $A^*_T(U) = A^*_T(X)/([S]_T)$
where $([S]_T)$ denotes the $T$-equivariant fundamental class
of $S$. Since $f$ has weight $12$, $[S]_T = 12t \in A^*_T(X) = Z[t]$.
Therefore, $A^*_T(U) = \Z[t]/12t$.

(b) The complement of $W$ in $X$ is the image
of the small diagonal $a = b = c$  under the degree $6$ map
$\A^3 \rightarrow \A^3$. The small diagonal has $T$-equivariant fundamental
class $4t^2$, the $T$ fundamental class of $X-W$ is $24t^2$. Hence
$A^*(\Mbar_{1,1}) = \Z[t]/24t^2$ as claimed. \endproof

\medskip

{\bf Remark.} From our computation we see $A^1({\cal M}_{1,1}) =
Pic_{fun}({\cal M}_{1,1}) = \Z/12$, a fact which was originally proved
by Mumford \cite{Mu}.
With an appropriate sign convention,
the class $t \in Z[t]/12t$ is just $c_1(L)$
where $L$ is the generator of $Pic_{fun}({\cal M}_{1,1})$ which
assigns to a family of elliptic
curves $X \stackrel{\pi} \rightarrow S$
the line bundle $\pi_*(\omega_{X/S})$ (the Hodge bundle).
Thus  monomial $at^n$ corresponds to the class which assigns to
a family ${\cal X} \rightarrow S$ of elliptic curves, the class
$c_1(L^{\otimes a})^n \cap [S] \in A_*(S)$.

Angelo Vistoli observed that this can be seen directly as follows:
The unipotent radical of $G$ acts freely on $U$ (or $W$) and
the quotient is the space of forms $y^2z = x^3 + \alpha xz^2 + \beta z^3$
with no double (or triple) roots. The torus
action is given by $t\cdot \alpha = t^{-4}\alpha$ and $t\cdot \beta
= t^{-6}\beta$, and the space is the total space of the $\G_m$
bundle over ${\cal M}_{1,1}$ (or $\overline{{\cal M}}_{1,1}$)
corresponding to the Hodge bundle. Thus, the
Chow ring of ${\cal M}_{1,1}$ (resp. $\overline{{\cal M}}_{1,1}$)
is generated by the first Chern class of the Hodge bundle.

\section{Some technical facts} \label{appendix}

\subsection{Intersection theory on algebraic spaces} \label{algspace}

Unfortunately, while most results about schemes generalize to algebraic
spaces, most references deal exclusively with schemes.  In particular,
this is the case for \cite{Fulton}, the basic reference for the
intersection theory used in this paper.  The purpose of this section
is to indicate very briefly how this theory generalizes to algebraic
spaces.

We recall from \cite{Knutson} the definition of algebraic spaces, and
basic facts about them.  If $X$ is a scheme, the functor $X^{\cdot} =
Hom(\cdot, X)$ from $(\mbox{Schemes})^{opp}$ to (Sets) is a sheaf in
either the Zariski or \'etale topologies.  With this as motivation, an
algebraic space is defined to be a functor $A: (\mbox{Schemes})^{opp}
\rightarrow (\mbox{Sets})$ such that:

(1) $A$ is a sheaf in the \'etale topology.

(2) (Local representability) There is a scheme $U$ and a
sheaf map $U^{\cdot} \rightarrow A$ such that for any scheme
$V$ with a map $V^{\cdot} \rightarrow A$, the fiber product
(of sheaves) $U^{\cdot} \times_{A} V^{\cdot}$ is represented
by a scheme, and the map
$U^{\cdot} \times_{A} V^{\cdot} \rightarrow V^{\cdot}$ is induced
by an \'etale surjective map of schemes.

Knutson also imposes a technical hypothesis of
quasi-separatedness, which states that the map $U^{\cdot}
\times_{A} U^{\cdot} \rightarrow U^{\cdot} \times U^{\cdot}$
is quasi-compact.

A morphism of algebraic spaces is a natural transformation of functors.

The map $X \mapsto X^{\cdot}$ is a fully faithful embedding of
(Schemes) into (Algebraic spaces).  We identify $X$ with $X^{\cdot}$
and henceforth use the same notation $X$ for both of these.  The
scheme $X$ is called a representable \'etale covering of $A$ (or an
\'etale atlas for $A$); it can be chosen to be a disjoint union of affine
schemes.  Thus, just as a scheme has a Zariski covering by affine
schemes, an algebraic space has an \'etale covering by affine schemes.

There are several ways to think of algebraic spaces in relation to
schemes.  One way is to think of a (normal) algebraic space as
a quotient of a scheme by a finite group (see \cite{Kollar}).  Another is
to think of an algebraic space as a quotient of a scheme by an \'etale
equivalence relation.  More precisely, in the setting of the above
definition, let $R$ denote the scheme $U \times_{A} U$; then $A$ is a
categorical quotient of $R \rightarrow U \times U$
\cite[II.1.3]{Knutson}.  Finally, any algebraic space has an
open dense subset which is isomorphic to a scheme
\cite[II.6.7]{Knutson}.

A key fact of algebraic spaces that we use
is:
\begin{prop} \label{l.algspacequotient}
Let $X$ be an algebraic space
with a set theoretically free action of $G$. Then a
quotient
$X \rightarrow Y$ exists in the category of algebraic spaces.  Moreover
this quotient is a principal bundle.
\end{prop}
Sketch of proof. Consider the functor $Y=[X/G]$ whose sections
over $B$ are principal $G$-bundles $E \rightarrow B$ together
with an equivariant map $E \rightarrow X$.
Since principal bundles can be constructed locally in the \'etale
topology, and $G$ acts
without stabilizers, one can check that
$Y$ is a sheaf in the \'etale topology (cf. \cite[Example
4.8]{D-M}).

To show that $Y$ is an algebraic space, we must construct an \'etale
atlas for $Y$. This follows from
\cite[Theorem 4.21]{D-M}. We sketch
a proof assuming that $X$ is normal.
Consider the surjective map of \'etale sheaves $X
\rightarrow Y$. It suffices to show that
every closed point $x \in X$ is contained in a locally closed subscheme $Z
\subset X$ such that $Z \rightarrow Y$ is \'etale.
Since $X$ is an algebraic space and
we are working locally in the \'etale topology
we may assume $X$ is a scheme.

Let $Gx \simeq G$
be the the $G$ orbit of $X$. Then $Gx$ is the fiber of $X \rightarrow
Y$ containing $x$. Let $Z$ be a locally closed subscheme of $X$
defined by lifts to ${\cal O}_X$ of the local equations for $x \in
Gx$ so that $\mbox{dim }Z = \mbox{dim }X - \mbox{dim }G$
and the scheme theoretic intersection $Z \cap Gx$ is $x \in Gx$.
Since $G$ is smooth, the point $x \in G$ is cut
out by a regular sequence of length $g$ where
$g = \mbox{dim }G$.

Now consider the equivariant map (where $G$ acts on $G \times
Z$ by $g(g_1,z) = (gg_1,z)$)
$G \times Z \stackrel{\psi} \rightarrow  X$ given by $(g,z) \mapsto gz$.
This map is the restriction to $G \times Z$ of
the action map $G \times X \stackrel{\Psi} \rightarrow X$. Since
$G$ acts without stabilizers, the fiber of $\Psi$
over $x$ is $\{(g^{-1},gx)\} \simeq Gx$.
The fiber of $\psi$ over $x$
is then
$(G \times Z \cap Gx)$. By construction of $Z$,
$Z \cap Gx \simeq x$. Thus, the scheme-theoretic fiber of $\psi$
over $x$ is $(1,x)$, so $\psi$ is unramified $x$.

Thus we have a surjection of complete local rings
${\cal O}_{x,X} \rightarrow {\cal O}_{(1,x),G \times Z} \rightarrow 0.$
Since $X$ and $G \times Z$ have the same dimension, and we assume
$X$ is normal, the map is also an injection since
$\widehat{\cal O}_{x,X}$ is integral domain. Therefore
$\psi$ is \'etale at $x$.
Since $G$ acts by automorphisms,
the open neighborhood of $(1,x)$ where $\psi$ is \'etale
is $G$-invariant. The group acts transitively on itself
so any invariant neighborhood of $(1,x)$ is of the
form $G \times U$ where $x \in U \subset Z$.
Since $G \times U \rightarrow X$ is \'etale
the map $U \rightarrow Y$ is \'etale.
Thus $Y =[X/G]$ has an \'etale cover by schemes and
is therefore an algebraic space.

Finally, $X
\rightarrow Y$ is a principal bundle, since $Y = [X/G]$ and
(tautologically) we have $X \times_{[X/G]} X = X \times G$.
\endproof \medskip

Using representable \'etale coverings, one can extend basic
definitions about schemes to algebraic spaces.  Much of this is done
in \cite{Knutson}, where more complete definitions and details can be
found.  Here are some examples.  Any sheaf on the category of schemes
(e.g. ${\cal O}$) extends uniquely to a sheaf on the category of
algebraic spaces: if $U$ is an affine \'etale covering of the space
$A$, and $R$ is as above, then ${\cal O}(A) = \mbox{Ker}( {\cal O}(U)
\rightarrow {\cal O}(R))$.  Likewise, a property $P$ of schemes is
called stable if given an \'etale covering $\{X_i \rightarrow X \}$,
$X$ has $P$ if and only if $X_i$ has $P$.  Any stable property of
schemes extends to a property of algebraic spaces by defining it in
terms of representable \'etale coverings.  Thus, one can speak of
algebraic spaces which are normal, smooth, reduced, $n$-dimensional,
etc.  Similarly, if $P$ is a stable property of maps of schemes such
that $P$ either (a) is local on the domain, or (b) satisfies effective
descent, then $P$ extends to a property of maps of algebraic spaces.
For example, one can speak of maps of algebraic spaces which are (a)
faithfully flat, flat, \'etale, universally open, etc., or (b) open
immersions, closed immersions, affine or quasi-affine morphisms, etc.

Likewise, again using representable \'etale coverings, one can extend
facts and constructions about schemes, e.g. $\mbox{Proj}$, fiber
products, divisors, etc., to algebraic spaces; again much of this is
done in \cite{Knutson}.

The definition of Chow groups of schemes given in \cite{Fulton}
generalizes immediately to algebraic spaces.  (A similar definition
was given for stacks in \cite{Gi}.)  If $X$ is an algebraic space,
define the group of $k$-cycles $Z_k(X)$ to be the free abelian group
generated by integral subspaces of dimension $k$.  To define rational
equivalence, first note that if $X$ is an integral algebraic space,
then the group of rational functions on $X$ is defined.  (Indeed, by
the above remarks, $X$ has an open dense subspace $X^0$ which is a
variety, and the rational functions on $X$ are the same as those on
$X^0$.)  If $Y$ is an integral subspace of $X$ of codimension $1$ and
$f$ is a rational function on $X$, then the order of vanishing of $f$
along $Y$, denoted $\mbox{ord}_Y(f)$, can be defined by taking an
\'etale map $\phi: U \rightarrow X$, where $U$ is a variety and
$\phi(U)$ has nonempty intersection with $Y$, and setting
$\mbox{ord}_Y(f) = \mbox{ord}_{\phi^{-1}(Y)}(\phi^*f)$, where the
right hand side is the definition for schemes in \cite{Fulton}.  If
$W$ is a $k+1$-dimensional integral subspace of $X$ and $f$ is a
rational function on $W$, define $\mbox{div}(f) \in Z_k(X)$ to be
$\sum \mbox{ord}_Y(f) [Y]$ where the sum is over all codimension $1$
integral subspaces of $Y$.  Then, exactly as in \cite{Fulton}, define
$\mbox{Rat}_k(X)$ to be the subgroup of $Z_k(X)$ generated by all
$\mbox{div}(f)$, for $f$ and $W$ as above, and define the Chow groups
$A_k(X) = Z_k(X) / \mbox{Rat}_k(X)$.

The arguments of \cite[Chapters 1-6]{Fulton} can be carried over
almost unchanged to show that Chow groups of algebraic spaces have
the same functorial properties as Chow groups of schemes.  Many of
the facts about schemes needed to prove this are extended to algebraic
spaces in \cite{Knutson}.

As a illustration, we will discuss the construction of Gysin
homomorphisms for regular embeddings, which is
the central construction of the first six chapters
of \cite{Fulton}.  If $X \rightarrow Y$ is any
closed embedding of algebraic spaces then we can define the cone
$C_XY$ as for schemes, since the $\mbox{Spec}$ construction for
sheaves (in the \'etale topology) of ${\cal O}_Y$ algebras defines an
algebraic space over $Y$.  If $X \rightarrow Y$ is a regular embedding
of codimension $d$, then $C_XY = N_XY$ is a vector bundle of rank $d$.
We can then define a specialization homomorphism $Z_k(Y) \rightarrow
Z_k(C_XY)$ as in \cite[Section 5.2]{Fulton}.  The deformation to the
normal bundle construction of \cite[Section 5.1]{Fulton} goes through
unchanged, since the blow-up $Y \times \P^1$ along the subspace $X
\times \infty$ is defined in the category of algebraic spaces. (The
existence of blow-ups is a consequence of the $\mbox{Proj}$
construction for graded algebras over algebraic spaces.) Thus as in
Fulton, the specialization map passes to rational equivalence.

In particular if $X \rightarrow Y$ is a regular embedding of
codimension $d$, then the construction of \cite[Chapter 6]{Fulton}
goes through, and we obtain a (refined) Gysin homomorphism.  If $X$ is
a (separated) smooth algebraic space, then the diagonal map is a
regular embedding. Therefore, the integral Chow groups of $X$ have an
intersection product.

Note also that algebraic spaces have an operational Chow
ring with the same formal properties as that of \cite[Chapter
17]{Fulton}. This follows from the fact that the ordinary Chow groups
of algebraic spaces have the same functorial properties as Chow groups
of schemes. In particular, if $X$ is a smooth algebraic space
of dimension $n$, then $A^i(X) = A_{n-i}(X)$, with the map
defined as in \cite{Fulton}.

{\bf Remark.} For algebraic stacks which have automorphisms, the
diagonal is not a regular embedding in the sense we have defined. If
the stack is Deligne-Mumford, then the diagonal is a local embedding
(i.e. unramified). For such morphisms, Vistoli constructed a Gysin
pullback with $\Q$ coefficients on the Chow groups of integral
substacks.  To obtain a good intersection theory with $\Z$ coefficients
on arbitrary algebraic stacks, a different definition of Chow groups
is required. For quotient stacks, equivariant Chow groups give a good
definition.  This point is discussed in Section \ref{itstacks}.

\subsection{Actions of group schemes over a Dedekind domain} \label{mixed}
Let $R$ be a Dedekind domain, and set $S = Spec(R)$.
Fulton's intersection theory remains valid for schemes and
thus algebraic spaces defined over $S$
(\cite[Section 20.2]{Fulton}). Thus, the
equivariant theory will work for actions of smooth affine
group schemes over $S$, provided we can find finite-dimensional
representations of $G/S$ where $G$ acts generically freely.
The following lemma shows that this can always be done if
the fibers of $G/S$ are connected.

\begin{lemma} \label{ded}
Let $G/S$ be a smooth affine group scheme defined over $S =\mbox{Spec }R$,
where $R$ is Dedekind domain. Then there exists a finitely generated
projective $S$-module $E$ , such that $G/S$ acts freely on an open
set $U \subset
E$ whose complement has arbitrarily high codimension.
\end{lemma}
Proof.
By \cite[Lemma 1]{Seshadri1} the coordinate ring $R(G)$
is a projective $R$-module with a $G$ action.
The group $G$ embeds into a finitely generated
$R(G)$ submodule $F$. Since $R$ is a Dedekind domain
$F$ is also projective. By \cite[Proposition 3]{Seshadri}
$F$ is contained in an invariant finitely generated submodule
$E$ (which is also projective).
$G$ acts freely on itself it acts freely
on an open of each fiber of $E/S$.
Replacing $E$ by $E \times_S
\ldots \times_S E$ we obtain a representation on which $G$ acts freely
on an open set $U \subset E$, such that $E-U$ has arbitrarily high
codimension.
\endproof \medskip

Thus if  $X/S$ is an algebraic space over $S$, we can construct
a mixed space $X_G = X \times_G U$, where $U$ is as in
the lemma.  We then define the
$i$-th equivariant Chow group as $A_{i+l-g}(X_G)$ where $l = \mbox{dim
}(U/S)$ and $g = \mbox{dim }(G/S)$.  Since most of the results of
intersection theory hold for algebraic spaces over a Dedekind
domain, most of the results on equivariant Chow
groups also hold, including the following:

(1) The functorial properties with respect to
proper, flat and l.c.i maps hold.

(2) If $X/S$ is smooth, there
is an intersection product on $A_*^G(X)$ for $X/S$ smooth.

(3) If $G$ acts freely on $X$ with quotient $X \rightarrow Y$
then $A_{*+g}^G(X) = A_*(Y)$.

(4) If $G/S$ acts with (locally) properly on $X/S$, then the theorem of
\cite{Kollar}, \cite{KM} implies that a quotient $X \rightarrow Y$
exists as an algebraic space over $S$. The results
of Section \ref{qint} (for ordinary Chow groups) generalize,
and
$A_{*+g}^G(X)_{\Q} = A_*(Y)_{\Q}$.

{\bf Remark.}
Facts (3) and (4) imply that any moduli space over $\mbox{Spec \Z}$
which is the quotient of a smooth algebraic space by a proper
action has a rational Chow ring.

\subsection{Some facts about group actions and quotients} \label{whenascheme}
Here we collect some useful results about actions of algebraic groups.

\begin{lemma} \label{l.free}
Suppose that $G$ acts properly on an an algebraic space
$X$. If the stabilizers are trivial, then the action is
free.
\end{lemma}
Proof. We must show that the action map $G \times X \rightarrow X \times X$
is a closed embedding. The properness of the action implies
that this map is proper and quasi-finite, hence finite. Since
the stabilizers are trivial, the map is unramified so it is an embedding
in a neighborhood of every point of $G \times X$. Finally, the map
is set theoretically injective, hence an embedding. \endproof \medskip

\begin{lemma} \label{q.exist} (\cite{E-G})
Let $G$ be an algebraic group. For
any $i > 0$, there is a representation $V$ of $G$ and an open set
$U \subset V$ such that $V-U$ has codimension more than $i$
and such that a principal bundle quotient $U \rightarrow U/G$
exists in the category of schemes.
\end{lemma}
Proof. Embed $G$ into $GL(n)$ for some $n$. Assume that
$V$ is a representation
of $GL(n)$ and $U \subset V$ is an open set such that a principal
bundle quotient $U \rightarrow U/GL(n)$ exists. Since $GL(n)$
is special, this principal bundle is locally trivial in the Zariski
topology. Thus $U$ is locally isomorphic to $W \times GL(n)$ for
some open $W \subset U/GL(n)$. A quotient $U/G$ can be constructed
by patching the quotients $W \times GL(n) \rightarrow W \times (GL(n)/G)$.
(It is well-known that a quotient $GL(n)/G$ exists \cite{Borel}.)

We have thus reduced to the case $G=GL(n)$. Let $V$ be
the vector space of $n \times p$ matrices with $p > i+ n$,
and let $U \subset V$ be the open set of matrices of maximal
rank. Then $V - U$ has codimension $p - n + 1$ and
the quotient $U/G$ is the Grassmannian $Gr(n,p)$.
\endproof \medskip

\medskip

The following proposition gives conditions under which the mixed space
$X_G$ is a scheme.  Recall that a group is special if every principal
bundle is locally trivial in the Zariski topology.  The groups
$GL(n)$, $SL(n)$, $Sp(2n)$, as well as solvable groups, are special;
$PGL(n)$ and $SO(n)$, as well as finite groups, are not \cite{Sem-Chev}.

\begin{prop} \label{inap}
Let $G$ be an algebraic group, let $U$ be a scheme on which $G$ acts
freely, and suppose that a principal bundle quotient $U \rightarrow
U/G$ exists.  Let $X$ be a scheme with a $G$-action.
Assume that one of the following conditions holds:

(1) $X$ is (quasi)-projective with a linearized $G$-action, or

(2) $G$ is connected and $X$ is equivariantly embedded as a closed
subscheme of a normal variety, or

(3) $G$ is special.

\noindent Then a principal bundle quotient
$X \times U \rightarrow X \times^G U$ exists in the category of
schemes.
\end{prop}
Proof. If $X$ is quasi-projective with a linearized action, then there
is an equivariant line bundle on $X \times U$ which is relatively
ample for the projection $X \times U \rightarrow U$. By \cite[Prop. 7.1]{GIT}
a principal bundle quotient $X \times^G U$ exists.

Now suppose that $X$ is normal and $G$ is connected.
By Sumihiro's theorem \cite{Sumihiro}, $X$ can be covered
by invariant quasi-projective open sets which have a linearized
$G$-action. Thus, by \cite[Prop. 7.1]{GIT} we can construct
a quotient $X_G = X \times^{G} U$ by patching the quotients of
the quasi-projective open sets in the cover.

If $X$ equivariantly embeds in a normal variety $Y$, then by the above
paragraph a principal bundle quotient $Y \times U \rightarrow Y
\times^G U$ exists. Since $G$ is affine, the quotient map is affine,
and $Y \times U$ can be covered by affine invariant open sets.  Since
$X \times U$ is an invariant closed subscheme of $Y \times U$, $X
\times U$ can also be covered by invariant affines. A quotient $X
\times^G U$ can then be constructed by patching the quotients of the
invariant affines.

Finally, if $G$ is special, then $U \rightarrow U/G$ is a
locally trivial bundle in the Zariski topology. Thus
$U = \bigcup\{U_\alpha\}$ where
$\phi_{\alpha}:U_\alpha \simeq G \times W_\alpha$
for some open $W_\alpha \subset U/G$. Then $\psi_{\alpha}:
X \times U_{\alpha} \rightarrow X \times W_{\alpha}$ is a quotient,
where $\psi_{\alpha}$ is defined by the formula
$(x,w,g) \mapsto (g^{-1}x,w)$.
(Here we assume that $G$ acts on the left on both factors
of $X \times G$.)  \endproof

\end{document}